\input amstex.tex
\input amsppt.sty
\magnification1200
\nologo
\voffset-.5cm
\hoffset.5cm


\def\ad{\operatorname{ad}}
\def\al{\alpha}

\def\be{\beta}
\def\bx{{\boxed{\phantom{\square}}\kern-.4pt}}

\def\Cb{\bar C}
\def\CC{{\Bbb C}}
\def\CMN{\CC^{\tss m}\ot\CC^N}
\def\crc{{\raise.24ex\hbox{$\sssize\kern.1em\circ\kern.1em$}}}

\def\d{\partial}
\def\de{\delta}
\def\De{\Delta}
\def\det{\operatorname{det}}

\def\End{\operatorname{End}\hskip1pt}
\def\EndCM{\End(\CC^M)}
\def\EndCN{\End(\CC^N)}
\def\enddemos{{\ $\square$\enddemo}}
\def\ep{\varepsilon}
\def\Esp{E^{\,\prime}}
\def\Et{\tilde E}

\def\Fsp{F^{\,\prime}}
\def\Ft{\tilde F}

\def\g{\frak{g}}
\def\ga{\gamma}
\def\gap{\gamma^{\tss\prime}}
\def\glM{\frak{gl}_M}
\def\glN{\frak{gl}_N}

\def\h{\frak{h}}
\def\Hf{\operatorname{Hf}}

\def\I{\Cal I}
\def\id{\operatorname{id}}
\def\io{\iota}
\def\Is{{I^\ast}}

\def\la{\lambda}
\def\La{\Lambda}

\def\lap{{\la^\prime}}
\def\las{{\la^\ast}}

\def\om{\omega}
\def\Om{\Omega}
\def\ot{\otimes}

\def\P{\Cal P}
\def\PD{\Cal P\Cal D}
\def\per{\operatorname{per}}
\def\Pf{\operatorname{Pf}}
\def\phi{\varphi}
\def\primes{{\tss\prime}}

\def\RR{{\Bbb R}}
\def\Rt{\tilde R}

\def\S{\operatorname{S}}
\def\sgn{\operatorname{sgn}\tss}
\def\si{\sigma}
\def\so{\frak{so}}
\def\soM{\so_M}
\def\soN{\so_N}
\def\sp{\frak{sp}}
\def\spM{\sp_M}
\def\spN{\sp_N}

\def\tr{\operatorname{tr}}
\def\tss{\hskip1pt}

\def\U{\operatorname{U}}
\def\Ug{\U(\g)}

\def\UglN{\U(\glN)}
\def\UsoM{\U(\soM)}
\def\UsoN{\U(\soN)}
\def\UspM{\U(\spM)}
\def\UspN{\U(\spN)}

\def\Z{\operatorname{Z}}
\def\Zg{\Z(\g)}

\def\ZglN{\Z(\glN)}
\def\ZsoM{\Z\tss(\tss\soM)}
\def\ZsoN{\Z\tss(\tss\soN)}
\def\ZspM{\Z\tss(\tss\spM)}
\def\ZspN{\Z\tss(\tss\spN)}
\def\ZZ{{\Bbb Z}}


\expandafter\ifx\csname bethe.def\endcsname\relax \else\endinput\fi
\expandafter\edef\csname bethe.def\endcsname{%
 \catcode`\noexpand\@=\the\catcode`\@\space}
\catcode`\@=11

\mathsurround 1.6pt
\font\Bbf=cmbx12 scaled 1200

\def\hcor#1{\advance\hoffset by #1}
\def\vcor#1{\advance\voffset by #1}
\let\bls\baselineskip  \let\ignore\ignorespaces
\def\vsk#1>{\vskip#1\bls} \let\adv\advance 
\def\vv#1>{\vadjust{\vsk#1>}\ignore} \def\vvv#1>{\vadjust{\vskip#1}\ignore}
\def\vvn#1>{\vadjust{\nobreak\vsk#1>\nobreak}\ignore}
\def\vvvn#1>{\vadjust{\nobreak\vskip#1\nobreak}\ignore}
\def\setnormalbls{\edef\normalbls{\bls\the\bls}}
\def\setmaths{\edef\maths{\mathsurround\the\mathsurround}}

 \let\nt\noindent 
\def\nn#1>{\noalign{\vskip #1pt}} \def\NN#1>{\openup#1pt}
 
\let\Sum\sum \def\sum{\Sum\limits} 
\let\Prod\prod \def\prod{\Prod\limits} \let\Int\int \def\int{\Int\limits}

\let\=\m@th \def\&{.\kern.1em} \def\>{\!\;} \def\:{\!\!\;}

\ifx\plainfootnote\undefined \let\plainfootnote\footnote \fi
\expandafter\ifx\csname amsppt.sty\endcsname\relax
 
\else \fi

\newbox\s@ctb@x
\def\s@ct#1 #2\par{\removelastskip\vsk>
 \vtop{\bf\setbox\s@ctb@x\hbox{#1} \parindent\wd\s@ctb@x
 \ifdim\parindent>0pt\adv\parindent.5em\fi\item{#1}#2\strut}%
 \nointerlineskip\nobreak\vtop{\strut}\nobreak\vsk-.4>\nobreak}

\newbox\t@stb@x
\def\gadv{\global\advance} \def\gad#1{\gadv#1 1} 
\def\l@b@l#1#2{\def\n@@{\csname #2no\endcsname}%
 \if *#1\gad\n@@ \expandafter\xdef\csname @#1@#2@\endcsname{\the\Sno.\the\n@@}%
 \else\expandafter\ifx\csname @#1@#2@\endcsname\relax\gad\n@@
 \expandafter\xdef\csname @#1@#2@\endcsname{\the\Sno.\the\n@@}\fi\fi}
\def\l@bel#1#2{\l@b@l{#1}{#2}\?#1@#2?}
\def\?#1?{\csname @#1@\endcsname}
\def\[#1]{\def\n@xt@{\ifx\t@st *\def\n@xt####1{{\setbox\t@stb@x\hbox{\?#1@F?}%
 \ifnum\wd\t@stb@x=0 {\bf???}\else\?#1@F?\fi}}\else
 \def\n@xt{{\setbox\t@stb@x\hbox{\?#1@L?}\ifnum\wd\t@stb@x=0 {\bf???}\else
 \?#1@L?\fi}}\fi\n@xt}\futurelet\t@st\n@xt@}
\def\(#1){{\rm\setbox\t@stb@x\hbox{\?#1@F?}\ifnum\wd\t@stb@x=0 ({\bf???})\else
 (\?#1@F?)\fi}}
\def\dff{\expandafter\d@f} \def\d@f{\expandafter\def}
\def\edff{\expandafter\ed@f} \def\ed@f{\expandafter\edef}

\newcount\Sno \newcount\Lno \newcount\Fno
\def\Section#1{\gadno\Fno=0\Lno=0\s@ct{\the\Sno.} {#1}\par} \let\Sect\Section
\def\section#1{\gad\Sno\Fno=0\Lno=0\s@ct{} {#1}\par} \let\sect\section
\def\l@F#1{\l@bel{#1}F} \def\<#1>{\l@b@l{#1}F} \def\l@L#1{\l@bel{#1}L}
\def\Tag#1{\tag\l@F{#1}} \def\Tagg#1{\tag"\llap{\rm(\l@F{#1})}"}
\def\Th#1{Theorem \l@L{#1}} \def\Lm#1{Lemma \l@L{#1}}
\def\Prop#1{Proposition \l@L{#1}}
\def\Cr#1{Corollary \l@L{#1}} \def\Cj#1{Conjecture \l@L{#1}}
 
\def\Proof#1.{\demo{\it Proof #1}}

 \def\setparindent{\edef\Parindent{\the\parindent}}
\def\Appendix{\Sno=64\let\p@r@\z@ 
 \def\Section##1{\gad\Sno\Fno=0\Lno=0 \s@ct{} \hskip\p@r@ Appendix
\\the\Sno
  \if *##1\relax\else {.\enspace##1}\fi\par} \let\Sect\Section
 \def\section##1{\gad\Sno\Fno=0\Lno=0 \s@ct{} \hskip\p@r@ Appendix%
  \if *##1\relax\else {.\enspace##1}\fi\par} \let\sect\section
 \def\l@b@l##1##2{\def\n@@{\csname ##2no\endcsname}%
 \if *##1\gad\n@@
 \expandafter\xdef\csname @##1@##2@\endcsname{\char\the\Sno.\the\n@@}%
 \else\expandafter\ifx\csname @##1@##2@\endcsname\relax\gad\n@@
 \expandafter\xdef\csname @##1@##2@\endcsname{\char\the\Sno.\the\n@@}\fi\fi}}

 \let\x\times 
 
\let\le\leqslant \let\ge\geqslant
  \let\8\infty

\let\=\m@th  \def\_#1{_{\rlap{$\ssize#1$}}}

\def\lc{{,\ldots\hskip-0.1pt,\hskip1pt}}

\def\E(#1){\mathop{\hbox{\rm End}\,}(#1)} 

\def\1{^{-1}} \def\vst#1{{\lower2.1pt\hbox{$\bigr|_{#1}$}}}

\let\logo@\relax
\let\m@k@h@@d\makeheadline \let\m@k@f@@t\makefootline
\def\makeheadline{\ifnum\pageno=1\headline={\hfil}\fi\m@k@h@@d}
\def\makefootline{\ifnum\pageno=1\footline={\hfil}\fi\m@k@f@@t}


\line{\Bbf Capelli Identities for Classical Lie Algebras\hfill}
\bigskip\bigskip
\line{\bf Alexander Molev$^1$, Maxim Nazarov$^2$\hfill}
\bigskip
\line{$^1$School of Mathematics and Statistics,
University of Sydney,\hfill}
\line{\phantom{$^1$}Sydney NSW 2006, Australia\tss.
E-mail: alexm{\@}maths.usyd.edu.au\hfill}
\smallskip
\line{$^2$Department of Mathematics, University of York,
York YO1 5DD, England\tss.\hfill}
\line{\phantom{$^2$}E-mail: mln1{\@}york.ac.uk\hfill}
\bigskip\bigskip


\line{\bf Abstract\hfill}
\kern10pt
\nt
We extend the Capelli identity from the Lie algebra $\glN$
to the other classical Lie algebras $\soN$ and $\spN$.
We employ the theory of reductive \text{dual pairs due to Howe.}


\bigskip\bigskip
\line{\bf1.\ Introduction\hfill}
\section{\,}\kern-20pt

\nt
The Capelli identity is one of the most renowned results in classical
invariant theory. It gives a set of $N$ distinguished generators for the
centre of enveloping algebra $\UglN$ of the 
general linear Lie algebra.
For any positive integer $m$, consider the natural action of the complex
Lie group $GL_N$ in the ring $\P$ of polynomials on the vector
space $\CC^{\tss m}\ot\CC^N$. The resulting representation 
\text{$\UglN\to\PD$ by differential}
operators on $\CC^{\tss m}\ot\CC^N$ with polynomial coefficients
is faithful when $m\ge N$.

The image of the centre $\ZglN$ of the algebra $\UglN$ under this
representation is the subalgebra $\I\subset\PD$ of $GL_m\x GL_N$-invariant
differential operators.
Denote $l=\min(m,N)$. The algebra $\I$ has a distinguished set of generators
$\Om_1\lc\Om_l$ which are called the {\it Cayley operators} [W].
Let $x_{ai}$ with $a=1\lc m$ and $i=1\lc N$ be coordinates on
$\CC^{\tss m}\ot\CC^N$. Let $\d_{ai}$ be the corresponding partial
differentiations.
Then $\Om_k\in\I$ is defined for any index $k=1,2,\ldots$ as
$$
\frac1{k\tss!}\ \,\ \sum_{\si\in S_k}\,\,
\sum_{a_1\lc a_k}\sum_{i_1\lc i_k}
{\sgn(\si)}\cdot
x_{a_1i_1}\ldots x_{a_ki_k}\,
\d_{a_1i_{\si(1)}}\ldots\d_{a_ki_{\si(k)}}\,
\Tag{0.1}
$$
where $\sgn(\si)$ is the sign of the element $\si$ of the
symmetric group $S_k$. Here the indices $a_1\lc a_k$ and $i_1\lc i_k$
run through $1\lc m$ and $1\lc N$ respectively. The differential
operator \(0.1) can be rewritten as
$$
\sum_{a_1<\ldots<a_k}\,\sum_{i_1<\ldots<i_k}
\det\tss[\tss x_{a_pi_q}\tss]\tss
\det\tss[\tss\d_{a_pi_q}\tss]
$$
where the determinants are taken with respect to the indices
$p,q=1\lc k$. In particular, we have $\Om_k=0$ for  any $k>l$.

The Capelli identity [C] gives an explicit formula for a preimage in $\ZglN$
of the Cayley operator $\Om_k\in\I$. Let $E_{ij}$ be the
generators of the enveloping algebra $\UglN$ such that in the above
representation $\UglN\to\PD$
$$
E_{ij}\,\mapsto\,\sum_a\,x_{ai}\d_{aj}\,.
\Tag{4.0}
$$
Then $\Om_k$ is the image of the element from $\UglN$ given by
$$
\frac1{k\tss!}\ \,\ 
\sum_{\si\in S_k}\,\,
\sum_{i_1\lc i_k}\ 
{\sgn(\si)}\,\cdot\tss
\prod_s\ 
\bigl(\tss E_{i_si_{\si(s)}}+(s-1)\cdot\de_{i_si_{\si(s)}}\tss\bigr)
\Tag{0.3}
$$
where the index $s$ runs through $1\lc k$ and the factors corresponding to $s$
are arranged from left to right. We write $\de_{ij}$ 
for the Kronecker delta. Here~\(0.3) is the $k$-th {\it Capelli element}
of the centre $\ZglN$;
\text{see [\tss H\tss,\tss HU\tss] for further comments.}

The algebra $\I$ has another distinguished set of generators
$\Theta_1\lc\Theta_l$. Invariant differential operator
$\Theta_k$ is defined for any $k=1,2,\ldots$ as
$$
\frac1{k\tss!}\ \,\ 
\sum_{\si\in S_k}\,\,
\sum_{a_1\lc a_k}\sum_{i_1\lc i_k}
x_{a_1i_1}\ldots x_{a_ki_k}\,
\d_{a_1i_{\si(1)}}\ldots\d_{a_ki_{\si(k)}}\,.
\Tag{0.11}
$$
Let $d_1\lc d_m$ be the multiplicities of the numbers $1\lc m$
respectively in the sequence $(a_1\lc a_k)$. Further, let
$f_1\lc f_N$ be the multiplicities of the numbers $1\lc N$
in the sequence $(i_1\lc i_k)$. Then
the differential
operator \(0.11) can be rewritten as
$$
\sum_{a_1\le\ldots\le a_k}\,\sum_{i_1\le\ldots\le i_k}\ 
\frac
{\per\tss[\tss x_{a_pi_q}\tss]\tss\per\tss[\tss\d_{a_pi_q}\tss]}
{d_1\tss!\,\ldots\,d_m\tss!\, f_1\tss!\,\ldots\,f_N\tss!}
$$
where
$$
\per\tss[\tss x_{a_pi_q}\tss]=
\sum_{\si\in S_k}\,
x_{a_1i_{\si(1)}}\ldots\tss x_{a_ki_{\si(k)}}
$$
is the {permanent} of the matrix $[\tss x_{a_pi_q}\tss]$
with the indices $p,q=1\lc k$.

An explicit formula for a preimage in $\ZglN$
of the operator $\Theta_k$ has been found 
by the second author [\tss N1\tss].
Similarly to \(0.3), 
the differential operator $\Theta_k$ is the image of the element from 
$\UglN$
$$
\frac1{k\tss!}\ \,\ 
\sum_{\si\in S_k}\,\,
\sum_{i_1\lc i_k}\ 
\tss
\prod_s\ 
\bigl(\tss E_{i_si_{\si(s)}}-(s-1)\cdot\de_{i_si_{\si(s)}}\tss\bigr)
\Tag{0.31}
$$
where the factors 
corresponding to the index 
$s=1\lc k$ are again arranged from left to right.
One can check that the element \(0.31) belongs to $\ZglN$; cf. [N1].
See also the remarks subsequent to the proof of Proposition 3.9
in the present article. 

Take an irreducible polynomial $\glN$-module of highest weight
\text{$\nu=(\nu_1\lc\nu_N)$.}
We make the standard choice of a Borel subalgebra in $\glN$; see Section 6.
Firstly, observe that the eigenvalue of \(0.3) in this
module is the {\it elementary symmetric polynomial}
$$
e_k(\nu_1\lc\nu_N)=\sum_{i_1<\ldots<i_k}\nu_{i_1}\ldots\tss\nu_{i_k}
\nopagebreak
$$
plus lower terms. Eigenvalue of \(0.31) 
is up to lower terms
the {\it complete symmetric polynomial}
$$
h_k(\nu_1\lc\nu_N)=\sum_{i_1\le\ldots\le i_k}\nu_{i_1}\ldots\tss\nu_{i_k}\,.
\nopagebreak
$$
Secondly, by the definitions \(0.1) and \(0.11) the eigenvalues of
\(0.3) and \(0.31) in that module vanish if $\nu_1+\ldots+\nu_N<k$.
These two properties characterize each of the 
central elements \(0.3) and \(0.31)
of $\UglN$ uniquely; cf.
[\tss K1\tss,\tss K2\tss,\tss OO1\tss,\tss S\tss].

Now let $\g$ be the classical Lie algebra $\soN$ or $\spN$ over $\CC$.
The rank of the reductive Lie algebra $\g$ is $n=[\tss N/2\tss]$.  
Let $G$ be the corresponding classical Lie group $O_N$ or $Sp_N$ over $\CC$.
We will denote by $\Zg$ the ring of invariants in the
enveloping algebra $\Ug$ with respect to the adjoint action of $G$.
So the ring $\Zg$ coincides with the centre of $\Ug$ for
$\g=\so_{2n+1}\tss,\tss\sp_{2n}$ but is strictly contained in the centre
of $\Ug$ when $\g=\so_{2n}$. The principal aim of this article is to
construct in $\Zg$ the analogues of the elements \(0.3) and \(0.31)
from $\ZglN$; \text{cf. \![\tss OO2\tss].}

In Section 2 for each $k=1,2,\ldots$ we define the elements $C_k$ and
$D_k$ of $\Zg$ similarly to the above characterization of 
the elements \(0.3) and \(0.31) of $\ZglN$.
We will choose a Cartan decomposition
$\g=\frak{n}_{-}\oplus\h\oplus\frak{n}_{+}$ in that section.
Now let $\la=(\la_1\lc\la_n)$
be a dominant weight of the Cartan subalgebra~$\h$
with integral labels.
So $\la_1\lc\la_n$ satisfy the inequalities (2.14).
Let $V_\la$ be the irreducible $\g$-module of the highest weight $\la$.
Each of the elements $C_k$ and $D_k$ of $\Zg$ is defined by 
the following two conditions;
cf. [\tss OO2\tss,\tss S\tss].
First, the eigenvalues of $C_k$ and $D_k$ in $V_\la$ are respectively
$(-1)^k\tss e_k(\la_1^2\lc\la_n^2)$ and $h_k(\la_1^2\lc\la_n^2)$ plus lower
terms. Second, both of the elements $C_k$ and $D_k$ vanish in $V_\la$ if
$\la_1+\ldots+\la_{n-1}+|\la_n|<k$.
In particular, we have $C_k=0$ for $k>n$.

In Section 3 we give explicit formulas for the elements $C_k$ and $D_k$
in terms of the standard generators of the enveloping algebra $\Ug$; see 
Theorems 3.2 and 3.3 respectively. These formulas are motivated by the
results of [\tss C1\tss,\tss S2\tss] and [\tss NO\tss].
We also employ the technique developed in our joint article
with G.\,Olshanski [\tss MNO\tss].
In our proofs we follow an approach to the classical Capelli
identity from [\tss M2\tss].

We will take the Lie group $G$ as a subgroup in $GL_N$.
This embedding is also fixed in Section 2. Let us restrict
the representation \(4.0) to the subalgebra $\Ug$ in $\UglN$.
In Section 4 we evaluate the images of the elements $C_1\lc C_n$
under this representation for $\g=\soN$; see Theorem 4.1. Further,
in Section 5 for $\g=\spN$ we evaluate the images 
under this representation of the elements
$D_1,D_2\,,\tss\ldots$; see Theorem 5.1. We regard these two results as
generalizations of the Capelli identity from $\glN$ to 
the Lie algebras $\soN$ and $\spN$ respectively.
A different generalization of
the Capelli identity has been
given in [\tss KS1\tss,\tss KS2\tss].

The algebra of invariants in $\PD$ with respect to
the natural action of the group $GL_N$ is generated
by the action of the Lie algebra $\frak{gl}_m$ in $\P$.
Consider the algebra of invariants in $\PD$ with respect to
the action of the subgroup $O_N\subset GL_N$.
This algebra is generated [\tss H\tss] by an action
in $\P$ of the Lie algebra $\spM$ with $M=2m$.
The images of $\ZsoN$ and $\ZspM$ in $\PD$
coincide. For each $k=1\lc m$ we express the image of
$C_k\in\ZspM$ in $\PD$ as
a linear combination of the images of $C_1\lc C_k\in\ZsoN$; see Theorem 4.4.
This result is based on the fact that the real Lie groups 
$O_N(\RR)$ and $Sp_M(\RR)$ form a reductive dual pair inside
the group $Sp_{MN}(\RR)$; cf. [\tss HU\tss,\tss KV\tss]. 
Along with Theorem~4.1 this result provides an explicit formula
for the image of $C_k\in\ZspM$ in $\PD$.   

Now consider the algebra of invariants in $\PD$ under the natural action
of the subgroup $Sp_N\subset GL_N$.
This algebra is generated by an action
in $\P$ of the Lie algebra $\soM$.
The images of $\ZspN$ and $\ZsoM$ in $\PD$ again
coincide. For each $k=1,2,\ldots$ we express the image of
$D_k\in\ZsoM$ in $\PD$ as
a linear combination of the images of $D_1\lc D_k\in\ZspN$; see Theorem 5.3.
It is based on the fact that the Lie groups 
$Sp_N(\RR)$ and $O_M^{\,\ast}(\RR)$ also form a reductive dual pair. 
\text{Along with} Theorem~5.1 it provides explicit formula
for the image of $D_k\in\ZspM$ in $\PD$.

In Section 6 we consider a certain generating function $C(u)$
for the elements $C_1\lc C_n\in\Zg$. This function 
is defined by (4.8) for $\g=\soN$ and by (5.10) for $\g=\sp_N$.
The main result of this section
is Theorem 6.2. It gives an explicit formula for the function $C(u)$
in terms of the standard generators of the algebra $\Ug$.
Moreover, $C(u)$ is invertible as a formal power series in $u\1$ and
the inverse series can be regarded as a generating function
for the elements $D_1,D_2,\tss\ldots\in\Zg$; see (4.12) 
if $\g=\soN$ and (5.11) if $\g=\spN$. Theorem 6.2 was obtained
in 1991 as a continuation of  [\tss O\tss] and served as a starting point
for the present work.

A $q$-analogue of the classical Capelli identity has been given in
[NUW1]. In the particular $m=1$ case $q$-analogues of
our Theorems 4.1 and 4.4
have been obtained in [\tss NUW2\tss].
It would be interesting to find $q$-analogues of all our results. 

We are grateful to A.\,Okounkov and G.\,Olshanski for kindly
informing us about
their results on the factorial Schur functions [\tss OO1\tss]
at an early stage of our work. We are also grateful to
R.\,Howe,
P.\,Jarvis,
T.\,Koornwinder,
L.\,Kov\'acs,
A.\,Lascoux,
B.\,Leclerc,
S.\,Okada,
C.\,Procesi
and J.-Y.\,Thibon
for helpful discussions and valuable remarks.
The second author has been supported by an 
EPSRC Advanced Research Fellowship and by EC TMR grant FMRX-CT97-0100.


\kern12pt
\line{\bf2.\ Two families of generators in $\Zg$\hfill}
\section{\,}\kern-20pt

\nt
In this section we will introduce two distinguished families
of generators in the ring $\Zg$ of $G$-invariants in the enveloping
algebra $\Ug$. Consider the vector space of symmetric polynomials in the
independent variables $z_1\lc z_n$ with complex coefficients
where $n$ is the rank of the simple
Lie algebra $\g$. We will make use of a special basis in this vector space
consisting of factorial Schur polynomials. Let us start with
recalling several facts about these polynomials from
\text{[\tss BL\tss,\tss M\tss,\tss OO1\tss].}

Fix an arbitrary sequence of complex numbers
$
a=(a_1,a_2,\dots).
$
Then for each $k=0,1,2,\ldots$ introduce the $k$-th
{\it generalized factorial power}
of a variable $z$
$$
(z\tss|\tss a)^k=(z-a_1)\cdots (z-a_k).
$$
Let $\mu=(\mu_1,\mu_2,\ldots)$ be a partition with length
$\ell(\mu)\le n$. Consider the function
$$
s_{\mu}(z_1\lc z_n\tss|\tss a)=
\frac
{\det\bigl[\tss(z_q|a)^{\mu_p+n-p}\tss\bigr]}
{\det\bigl[\tss(z_q|a)^{n-p}\tss\bigr]}
\Tag{2.1}
$$
where the determinants are taken with respect to the indices $p,q=1\lc n$.
This function is a symmetric polynomial in $z_1\lc z_n$ which is called
the {\it generalized factorial Schur polynomial\,}
[\tss M\tss,\tss Example I.3.20\tss].
Note that here the denominator
$$
\det\bigl[\tss(z_q|a)^{n-p}\tss\bigr]\,=\,\prod_{p<q}\,(z_p-z_q)
\Tag{2.15}
$$
is the Vandermonde determinant. Thus it does not depend
on the sequence $a$.

If $a=(0,0,\ldots)$ the function
$s_\mu(z_1\lc z_n\tss|\tss a)$ is the ordinary Schur polynomial
$s_\mu(z_1\lc z_n)$.
For the general sequence $a$ the definition \(2.1) implies that
$$
s_{\mu}(z_1\lc z_n\tss|\tss a)=s_\mu(z_1\lc z_n)+
\text{\rm lower degree terms}.
\nopagebreak
\Tag{2.2}
$$
Hence the polynomials $s_{\mu}(z_1\lc z_n\tss|\tss a)$
with $\ell(\mu)\le n$ form a linear basis
in the ring of symmetric polynomials in the variables $z_1\lc z_n$.

{}From now on we will suppose that the sequence $a$ is multiplicity-free,
that is $a_k\ne a_l$ for all $k\neq l$. 
We will use the following
characterisation of the polynomial $s_{\mu}(z_1\lc z_n\tss|\tss a)\,$; 
\text{cf. [\tss OO1\tss,\tss Theorem 3.3\tss].}
For any partition $\la=(\la_1,\la_2,\ldots)$ such that $l(\la)\le n$
introduce
the following $n$-tuple of elements from $a$:
$$
a_{\la}=(a_{\la_1+n}\lc a_{\la_n+1}).
$$
As usual, we will write $|\la|$ for the degree $\la_1+\la_2+\ldots$
of the partition $\la$.

\proclaim{Theorem 2.1}
Let $f(z_1\lc z_n)$ be a symmetric polynomial of degree not more than
$|\mu|$. Then the following three conditions are equivalent\tss:
\itemitem{\rm i)}
$f(z_1\lc z_n)$ equals $s_\mu(z_1\lc z_n\tss|\tss a)\,$
up to a 
scalar multiple\tss;
\itemitem{\rm ii)}
$f(a_\la)=0$ for any $\la$ such that $\la_k<\mu_k$ for at least one index
$k$;
\itemitem{\rm iii)}
$f(a_\la)=0$ for any $\la$ such that $|\la|<|\mu|$ while the leading
component of the polynomial $f(z_1\lc z_n)$ is $s_{\mu}(z_1\lc z_n)$
up to a scalar multiple.
\endproclaim

\demo{Proof}
Let us check first that the polynomial $s_{\mu}(z_1\lc z_n\tss|\tss a)\,$
satisfies (ii). For $(z_1\lc z_n)=a_\la$
the $(p,q)$-entry of the matrix
\text{in the numerator in \(2.1) is}
$$
(a_{\la_q+n-q+1}-a_1)\,\ldots\,(a_{\la_q+n-q+1}-a_{\mu_p+n-p}).
\Tag{1.8}
$$

Suppose that $\la_k<\mu_k$ for some index $k$. Then for $p\le k\le q$ we have
$$
1\le\la_q+n-q+1\le\la_k+n-k+1\le
\mu_k+n-k\le\mu_p+n-p\tss.
$$
Therefore every entry \(1.8) with $p\le k\le q$ is zero.
So the numerator in \(2.1) vanishes at $(z_1\lc z_n)=a_\la$.
The denominator \(2.15) does not
vanish at $a_\la$ since the sequence $a$ is multiplicity-free.
Hence $s_\mu(a_\la\tss|\tss a)=0$.

Let us also show that $s_\mu(a_\mu\tss|\tss a)\neq0$. Indeed, for $\la=\mu$
the $(p,q)$-entry \(1.8) equals
$$
(a_{\mu_q+n-q+1}-a_1)\,\ldots\,(a_{\mu_q+n-q+1}-a_{\mu_p+n-p})
$$
which is zero for $p<q$ and non-zero for $p=q$. So the square matrix in
the numerator of \(2.1) with $(z_1\lc z_n)=a_\mu$ is lower-triangular
with non-zero diagonal entries. So $s_\mu(a_\mu\tss|\tss a)\neq0$.

Note that the polynomial $s_{\mu}(z_1\lc z_n\tss|\tss a)$
satisfies the condition (iii); see \(2.2).

Let us now fix any symmetric polynomial $f(z_1\lc z_n)$ of degree not
more than $|\mu|$. Suppose that the condition (ii) 
is satisfied.
Let us represent this polynomial as a linear combination of
the polynomials $s_{\nu}(z_1\lc z_n\tss|\tss a)\tss$:
$$
f(x)=c_\mu\,s_{\mu}(z_1\lc z_n\tss|\tss a)+
\sum_{\nu\neq\mu} c_{\nu}\,s_{\nu}(z_1\lc z_n\tss|\tss a)\,;
\qquad c_\mu,c_\nu\in\CC.
\Tag{1.11}
$$
Here we assume that $|\nu|\le|\mu|$. Note that under this assumption
the condition $\nu\neq\mu$ is equivalent to existence of an index $k$
such that $\nu_k<\mu_k\tss$. Denote by $\Cal S$ the set of all possible
partitions $\nu$. Let the partition $\la$ also run
\text{through the set $\Cal S$.}
Put $(z_1\lc z_n)=a_\la$ in \(1.11). Then we obtain the system
of linear equations on the coefficients $c_\nu$
$$
\sum_{\nu\in\Cal S} c_{\nu}\,s_{\nu}(a_\la\tss|\tss a)=0\,;
\ \quad\la\in\Cal S.
$$
Equip the set $\Cal S$ with any linear ordering $\prec$ such that
$|\la|<|\nu|$ implies $\la\prec\nu$. By arranging
the rows and columns
of the matrix $[\tss s_{\nu}(a_{\la}|a)\tss]$ according to
the ordering $\prec$
we make this matrix triangular with non-zero
diagonal entries.
Hence $c_\nu=0$ for every $\nu\in\Cal S$.
So the polynomial $f(z_1\lc z_n)$ satisfies (i) by \(1.11).
The proof that the condition (iii) on $f(z_1\lc z_n)$
implies (i) is quite similar 
\enddemos

\nt
Let us now define the {\it factorial elementary}
and {\it complete} symmetric polynomials by
$$
\align
e_k(z_1\lc z_n\tss|\tss a)&=s_{(1^k)}(z_1\lc z_n\tss|\tss a),
\Tag{2.4}
\\
h_k(z_1\lc z_n\tss|\tss a)&=s_{(k)}(z_1\lc z_n\tss|\tss a)
\Tag{2.5}
\endalign
$$
respectively.
We also have the following explicit formulas for these polynomials.

\proclaim{Proposition 2.2}
We have the equalities
$$
\align
e_k(z_1\lc z_n\tss|\tss a)
\ &=\sum_{1\le p_1<\ldots <p_k\le n}(z_{p_1}-a_{p_1})
(z_{p_2}-a_{p_2-1})\cdots (z_{p_k}-a_{p_k-k+1})\,,
\\
h_k(z_1\lc z_n\tss|\tss a)
\ &=\sum_{1\le p_1\le\ldots \le p_k\le n}(z_{p_1}-a_{p_1})
(z_{p_2}-a_{p_2+1})\cdots (z_{p_k}-a_{p_k+k-1})\,.
\endalign
$$
\endproclaim

\demo{Proof}
We will make use of Theorem 2.1.
Denote respectively by $e_k^\primes(z_1\lc z_n\tss|\tss a)$ and 
$h_k^\primes(z_1\lc z_n\tss|\tss a)$
the polynomials on the right hand sides of the above relations.

Let us prove first that these polynomials are symmetric
in $z_1\lc z_n$. It suffices to demonstrate that they are invarant
under the transposition of $z_i$ and $z_{i+1}$ for each $i=1,\dots,n-1$.
But this reduces our task to the case $n=2$. Then we have 
$e_k^\primes(z_1,z_2\tss|\tss a)=0$ for $k>2$ while both
$e_1^\primes(z_1,z_2\tss|\tss a)$ and 
$e_2^\primes(z_1,z_2\tss|\tss a)$ are obviously symmetric.
The symmetry of the polynomial
$h_k^\primes(z_1,z_2\tss|\tss a)$
follows from the equality
$$
h_k^\primes(z_1,z_2\tss|\tss a)=
\frac{(z_1\tss|\tss a)^{k+1}-(z_2\tss|\tss a)^{k+1}}{z_1-z_2}
$$
which is easy to verify. Indeed, the product
$(z_1-z_2)\cdot h_k^\primes(z_1,z_2\tss|\tss a)$ equals
$$
\align
(z_1-z_2)\,\cdot\tss
\sum_{l=0}^k\,\,
(z_1-a_1)\ldots(z_1-a_l)\,(z_2-a_{l+2})\ldots(z_2-a_{k+1})
&=
\\
\sum_{l=0}^k\,\,
(z_1-a_1)\ldots(z_1-a_{l+1})\,(z_2-a_{l+2})\ldots(z_2-a_{k+1})
&\,\tss-
\\
\sum_{l=0}^k\,\,
(z_1-a_1)\ldots(z_1-a_l)\,(z_2-a_{l+1})\ldots(z_2-a_{k+1})
&\,\tss.
\endalign
$$
All summands at the right hand side of the above equality
cancel each other except
$$
(z_1-a_1)\ldots(z_1-a_{k+1})-(z_2-a_1)\ldots(z_2-a_{k+1})
=(z_1\tss|\tss a)^{k+1}-(z_2\tss|\tss a)^{k+1}\tss.
$$

Polynomials
$e_k(z_1\lc z_n|a),e_k^\primes(z_1\lc z_n|a)$ have the same
leading components
and so do the polynomials
$h_k(z_1\lc z_n\tss|\tss a),h_k^\primes(z_1\lc z_n\tss|\tss a)$.
Thus to complete the proof of Proposition 2.2 it remains to check
that $e_k^\primes(a_{\lambda}\tss|\tss a)=0$ and
$h_k^\primes(a_{\lambda}\tss|\tss a)=0$ for
$|\lambda|<k$.
By the symmetry 
of polynomials $e_k^\primes(z_1\lc z_n\tss|\tss a)$
and $h_k^\primes(z_1\lc z_n\tss|\tss a)$
$$
\align
e_k^\primes(z_1\lc z_n\tss|\tss a)
\ &=\sum_{1\le p_1<\ldots <p_k\le n}(z_{n-p_1+1}-a_{p_1})
\cdots (z_{n-p_k+1}-a_{p_k-k+1})\,,
\\
h_k^\primes(z_1\lc z_n\tss|\tss a)
\ &=\sum_{1\le p_1\le\ldots \le p_k\le n}(z_{n-p_1+1}-a_{p_1})
\cdots (z_{n-p_k+1}-a_{p_k+k-1})\,.
\endalign
$$
Note that $p_1\le n-k+1$ in the first of
these equalities. Thus if $\lambda_k=0$ then
$(a_{\lambda})_{n-p_1+1}=a_{p_1}$ which implies that
$e_k^\primes(a_{\lambda}\tss|\tss a)=0$.
Finally, consider the last formula for $h_k^\primes(z_1\lc z_n\tss|\tss a)$.
If $\lambda_1<k$ then 
$$
0\le \lambda_{n-p_1+1}\le \cdots \le \lambda_{n-p_k+1}\le k-1\tss.
\nopagebreak
$$
These inequalities imply that there exists $i\in\{1,\dots,k\}$ such that
$\lambda_{n-p_i+1}=i-1$. Then $(a_{\lambda})_{n-p_i+1}=a_{p_i+i-1}$
and $h_k^\primes(a_{\lambda}\tss|\tss a)=0$
\enddemos

\nt
In particular, $e_k(z_1\lc z_n\tss|\tss a)=0$ for $k>n$.
Consider the rational function in $t$ 
$$
X(t|a)=\frac{(t-z_1)\ldots (t-z_n)}{(t-a_1)\ldots (t-a_n)}\,.
$$
Next proposition shows that $X(t|a)$ and
$X(t|a)\1$ can be regarded as generating functions for the
families of symmetric polynomials \(2.4) and \(2.5).

\proclaim{Proposition 2.3}
We have the equalities of formal power series in $t\1$
$$
\align
1+\sum_{k=1}^n
\frac{(-1)^k e_k(z_1\lc z_n\tss|\tss a)}{(t-a_{n-k+1})\ldots (t-a_n)}
&=X(t|a)\,,
\Tag{1.6}
\\
1+\sum_{k=1}^{\infty}
\frac{h_k(z_1\lc z_n\tss|\tss a)}{(t-a_{n+1})\cdots (t-a_{n+k})}
&=X(t|a)\1\,.
\Tag{1.7}
\endalign
$$
\endproclaim

\demo{Proof}
We will employ the arguments from [\tss OO1\tss,\tss Theorem 12.1\tss].
Let us prove \(1.6) first. There is a unique decomposition
$$
\frac{(t-z_1)\ldots (t-z_n)}{(t-a_1)\ldots (t-a_n)}=
1+\sum_{k=1}^n\frac{c_k(z_1\lc z_n)}{(t-a_{n-k+1})\ldots (t-a_n)}
\Tag{2.33}
$$
where $c_k(z_1\lc z_n)$ is a certain symmetric polynomial in $z_1\lc z_n$
of degree $k$. Moreover, the leading component of this polynomial
coincides with that of
$(-1)^k e_k(z_1\lc z_n\tss|\tss a)$. Due to Theorem 2.1 it suffices
to prove that $c_k(a_\la)=0$ for any partition $\la$ with $|\la|<k$.
Let us substitute $(z_1\lc z_n)=a_\la$ with $\ell(\la)<k$ in \(2.33).
At the left hand side we get the function of $t$
$$
\prod_{p=1}^{k-1}\ \frac{t-a_{\la_p+n-p+1}}{t-a_{n-p+1}}\,.
\nopagebreak
$$
This rational function obviously has less than $k$ poles.
Therefore the coefficient $c_k(a_\la)$ at the right hand side vanishes.

The proof of the identity \(1.7) is similar. Consider the unique decomposition
$$
\frac{(t-a_1)\ldots (t-a_n)}{(t-z_1)\ldots (t-z_n)}=
1+\sum_{k=1}^{\infty}
\frac{d_k(z_1\lc z_n)}{(t-a_{n+1})\cdots (t-a_{n+k})}
\Tag{2.36}
$$
where $d_k(z_1\lc z_n)$ is a certain symmetric polynomial in $z_1\lc z_n$
of degree $k$. The leading component of this polynomial
coincides with that of
$h_k(z_1\lc z_n\tss|\tss a)$. It remains to prove that
$d_k(a_\la)=0$ for any partition $\la$ with $|\la|<k$.
Let us substitute $(z_1\lc z_n)=a_\la$ in \(2.36).  
At the left hand side we get a rational function in $t$
with simple poles only at $t=a_{\la_p+n-p+1}$ where $\la_p+n-p+1>n$,
that is where $\la_p\ge p$. The number of such indices $p$
cannot exceed $\la_1$. So the coefficient $d_k(\la)$ with $k>\la_1$
at the right hand side vanishes
\enddemos

\nt
We will use the symmetric polynomials \(2.4) and \(2.5) to define
some \text{distinguished} elements $C_k$ and $D_k$ in the ring of
invariants $\Zg$.
{}From now on we will let the indices $i\tss,j$ run through the set
$\{\tss-\tss n\lc\!-\!1,1\lc n\tss\}$ if $N=2n$ and the set 
$\{\tss-\tss n\lc\!-\!1,0,1\lc n\tss\}$ if $N=2n+1$.
Let $e_i$ form the standard basis in $\CC^N$. 
We will realize the group $G=O_N$ as the subgroup in $GL_N$ preserving
the symmetric bilinear form
$
\langle\,e_i\tss,e_j\,\rangle=\de_{i,-j}\,.
$
The group $G=Sp_N$
will be realized as the subgroup in $GL_N$ preserving the
alternating form
$
\langle\,e_i\tss,e_j\,\rangle=\de_{i,-j}\cdot\sgn i\,.
$

Put $\ep_{ij}=\sgn i\cdot\sgn j$ if $G=Sp_N$ and
\text{$\ep_{ij}=1$ if $G=O_N$.}
Let $E_{ij}\in\EndCN$ be the standard matrix units. We will also regard
$E_{ij}$ as basis elements of the Lie algebra $\glN$.
The Lie subalgebra $\g\in\glN$ is then spanned by the elements
$F_{ij}=E_{ij}-\ep_{ij}\cdot E_{-j,-i}$.
Now let us fix the Cartan decomposition
$\g=\frak{n}_{-}\oplus\h\oplus\frak{n}_{+}$ as follows.
The nilpotent subalgebras $\frak{n}_{+}$ and $\frak{n}_{-}$ are spanned
by the elements $F_{ij}$ where $i<j$ and $i>j$ respectively.
The Cartan subalgebra $\h$ is spanned by
the elements $F_{ii}$. Further, we will fix the basis
$\bigl(F_{-n,-n}\tss\lc\tss F_{-1,-1}\bigr)$ in $\h$.  
Any weight $\la=(\la_1\lc\la_n)$ of $\h$ will be taken with respect to this
basis. The half-sum of the positive roots of $\h$ is then
$$
\rho=(\ep+n-1,\ep+n-2\lc\ep)
\Tag{2.88888888}
$$
where $\ep=0\,,\,\frac12\,,\,1$ for $\g=\so_{2n}\,,\,\so_{2n+1}\,,\,\sp_{2n}$
respectively. Denote $\la_p+\rho_p=l_p$.

Recall that the Harish-Chandra homomorphism $\om:\Ug^\h\to\U(\h)$
is defined by
$$
\om\tss(X)-X\,\in\,
\Ug\,\frak{n}_{+}\cap\Ug^\h=
\frak{n}_{-}\!\Ug\cap\Ug^\h
$$
We will identify the algebra $\U(\h)=\S(\h)$ with the
algebra of polynomial functions on the space $\h^\ast$.
Then the restriction of $\om$ to the subalgebra $\Zg\subset\Ug^\h$
is an isomorphism onto the algebra of symmetric polynomials in the
variables $l_1^2\lc l_n^2$ [\tss D\tss,\tss Theorem 7.4.5\tss].
We will call this restriction the {\it Harish-Chandra isomorphism.}

Let us now make a choice of the sequence $a$ in the definition \(2.1).
Namely, set
$$
a=\bigl(\,\tss\ep^2,(\ep+1)^2,(\ep+2)^2,\tss\ldots\,\,\tss\bigr)\tss.
\Tag{2.8}
$$
Note that here $\{\tss a_1\lc a_n\}=\{\,\rho_1^2\lc\rho_n^2\,\}$.
Then define the elements $C_k$ and $D_k$ in $\Zg$ as the preimages
with respect to the Harish-Chandra isomorphism of the symmetric
polynomials $(-1)^k\tss e_k(l_1^2\lc l_n^2\tss|\tss a)$ and
$h_k(l_1^2\lc l_n^2\tss|\tss a)$ respectively. 
Note that $C_k=0$ for any $k>n$. Obviously, 
we have the following proposition.

\proclaim{Proposition 2.4}
Each of the families $C_1\lc C_n\!$ and $D_1\lc D_n\!$
\text{generates~$\Zg.\hskip-0.24pt$}
\endproclaim

\nt
In the next section we will give explicit formulas for the
elements $C_k,D_k\in\Zg$ for any $k$ in terms of the standard generators
$F_{ij}$ of the enveloping algebra $\Ug.$

We will use a corollary to Theorem 2.1.
Let $\la=(\la_1\lc\la_n)$ be any dominant weight of $\h$
with integral labels. So
$\la_1\lc\la_n$ satisfy the inequalities
$$
\alignat2
&\la_1\ge\ldots\ge\la_n\ge0&&\qquad\text{for}\quad
\g=\so_{2n+1}\,,\,\sp_{2n}\,;
\Tag{0.999}
\\
&\la_1\ge\ldots\ge\la_{n-1}\ge|\la_n|&&\qquad\text{for}\quad\g=\so_{2n}\,.
\endalignat
$$
In the case $\g=\so_{2n}$ we denote $\las=(\la_1\lc\la_{n-1},-\la_n)$.
Let $V_\la$ be the irreducible $\g$-module of the highest weight $\la$.
Note that for $\g=\so_{2n}$ the eigenvalues of any element from $\Zg$ in the
irreducible modules $V_\la$ and $V_\las$ coincide.
We will assume that $\la_n\ge0$ always, so that $\la$
will be always a partition with $\ell(\la)\le n$.

Consider the standard $\ZZ$-filtration on the enveloping algebra
$\Ug$. Let $Z$ be any element in the ring $\Zg$ of degree not more than $2k$.
Then here is the corollary.

\proclaim{Corollary 2.5}
\hskip-0.3pt Suppose that $Z$ 
vanishes in every module $V_\la$ with \text{$\ell(\la)<k$}
or in every module $V_\la$ with $\la_1<k$.
Then $Z$ equals respectively $C_k$ or $D_k$ up to a scalar multiple.
\endproclaim

\demo{Proof}
The eigenvalue of the element $Z\in\Zg$ is a certain symmetric polynomial
in $l_1^2\lc l_n^2$. Let us denote by $f$ this polynomial, its degree
does not exceed~$k$. By our choice \(2.8) of the sequence $a$ we have
$(\tss l_1^2\lc l_n^2)=a_\la$.
Thus $f(a_\la)=0$ for every $\la$ such that \text{$\ell(\la)<k$}
or $\la_1<k$. Theorem 2.1 now implies that  
$f(l_1^2\lc l_n^2)$ equals respectively $e_k(l_1^2\lc l_n^2\tss|\tss a)$
or $h_k(l_1^2\lc l_n^2\tss|\tss a)$ up to a scalar multiple for any $\la$.
But the latter two expressions are the eigenvalues in $V_\la$ of the
elements $(-1)^k\tss C_k$ and $D_k$ respectively. So we get the corollary
since any element of $\Zg$ is uniquely determined by 
its eigenvalues in the irreducible modules $V_\la$
\enddemos


\kern8pt
\line{\bf3.\ Explicit formulas for the elements $C_k$ and $D_k$\hfill}
\section{\,}\kern-20pt

\nt
In this section we employ the technique developed in
[\tss MNO\tss,\tss NO\tss] and \text{[\tss M1\tss,\tss M2\tss,\tss N2\tss].}

As in the previous section, let $E_{ij}\in\EndCN$
be the standard matrix units. But here we will regard
the elements $F_{ij}=E_{ij}-\ep_{ij}\cdot E_{-j,-i}$
as generators of the enveloping algebra $\Ug$. Let $u$ be a complex variable.
Introduce the element
$$
F=\sum_{ij}\,E_{ij}\ot F_{ji}\in\EndCN\ot\Ug
\Tag{3.00}
$$
and denote
$$
F(u)=F-u-\eta\,\in\,\EndCN\ot\Ug\tss[\tss u\tss]
\Tag{3.02}
$$
where we set $\eta=\frac12,-\frac12$ for $\g=\soN,\spN$ respectively.  
We will use this \text{notation} throughout the present section.
Consider the elements of the algebra $\EndCN^{\ot2}$
$$
P=\sum_{ij}\,E_{ij}\ot E_{ji}\,,
\ \quad
Q=\sum_{ij}\,\ep_{ij}\cdot E_{ij}\ot E_{-i,-j}\,.
\nopagebreak
\Tag{3.025}
$$
The element $P$ corresponds to the exchange operator
$e_i\ot e_j\mapsto e_j\ot e_i$ in $(\CC^N)^{\ot2}$.
The element $Q$ is obtained from $P$ by applying
to either of the tensor factors in $\EndCN^{\ot2}$ 
the transposition with respect to the bilinear form $\langle\,\,,\,\rangle$.
Observe that
$$
P\tss Q=Q\hskip.5pt P=
\cases
\phantom{-\,}Q&\text{if\ \ $\g=\soN$,}
\\
-\,Q&\text{if\ \ $\g=\spN$.}
\endcases
\Tag{3.01}
$$
Denote
$$
R(u,v)=1-\frac{P}{u-v}\,,
\qquad
\Rt(u,v)=1+\frac{Q}{u+v} 
$$
in $\!\EndCN^{\ot2}(u,v)$.
\hskip-0.332pt\!The first of these two functions is the
\text{\it rational Yang \hskip-0.5pt$R$\hskip-0.64pt-matrix\,}.

We will denote by $\io_p$ the embedding of the algebra $\EndCN$
into a finite tensor product $\EndCN^{\ot m}$ as the $p$-th tensor factor:
$$
\io_p(X)=1^{\ot\tss(p-1)}\ot X\ot1^{\ot\tss(m-p)}\tss;
\qquad
p=1\lc m\tss.
\nopagebreak
$$
Then put
$$
F_p(u)=(\io_p\ot\id)\,\bigl(F(u)\bigr)\in
\EndCN^{\ot m}\ot\Ug\,[\tss u\tss]\tss.
\Tag{3.0003}
$$
Direct calculation yields the next proposition; 
for details see \text{[\tss MNO\tss,\tss Section 3.11\tss].}
In this proposition the elements $R(u,v)$ and $\Rt(u,v)$
of $\EndCN^{\ot2}(u,v)$ are identified
with the elements $R(u,v)\ot1$ and $\Rt(u,v)\ot1$
of $\EndCN^{\ot2}\ot\Ug\,(u,v)$ respectively.
We will keep using this convention in the present section for simplicity.

\proclaim{Proposition 3.1}
We have the relation in the algebra $\EndCN^{\ot2}\ot\Ug\,(u,v)$
$$
R(u,v)\,F_1(u)\,\Rt(u,v)\,F_2(v)=F_2(v)\,\Rt(u,v)\,F_1(u)\,R(u,v)\,.
\Tag{3.111}
$$
\endproclaim

\nt
We will also use various embeddings of the algebra $\EndCN^{\ot2}$
into $\EndCN^{\ot m}$ with $m\ge2$. For any $1\le p<q\le m$ and   
$Y\in\EndCN^{\ot2}$ we will denote
$$
Y_{pq}=(\io_p\ot\io_q)\,(Y)\tss\in\tss\EndCN^{\ot m}\tss.
\Tag{3.002255}
$$
In this section the number $m$ of the tensor factors will also vary,
but $\mu$ will be always one of the two partitions $(1^m),(m)$.
Let $A_{m}$ and $B_m$ be the elements of $\EndCN^{\ot m}$
corresponding
to antisymmetrization and symmetrization in
the tensor product $(\CC^N)^{\ot m}$. They are normalized so that
$A_{m}^{\tss2}=A_{m}$ and $B_{m}^{\tss2}=B_{m}$.
Introduce the rational function $F_\mu(u)$
valued in $\EndCN^{\ot m}\ot\Ug$ as follows:
$$
\align
F_{(1^m)}(u)
&=(A_m\ot1)\,\cdot\,\prod_{q=1}^m\ 
\biggl(\!
\Bigl(\tss1+\frac{Q_{1q}+\ldots+Q_{q-1,q}}{2u-q+1}\tss\Bigr)\ot1
\tss\!\biggr)
\,F_q(u-q+1)
\tss,
\\
F_{(m)}(u)
&=(B_m\ot1)\,\cdot\,\prod_{q=1}^m\ 
\biggl(\!
\Bigr(\tss1+\frac{Q_{1q}+\ldots+Q_{q-1,q}}{2u+q-1}\tss\Bigr)\ot1
\biggr)
\,F_q(u+q-1)
\tss;
\endalign
$$
we will always arrange the (non-commuting)
factors corresponding to the indices $q=1\lc m$
from left to right.

Let us now choose the {\it classical point} \text{$u_\mu\in\CC$ as}
$$
u_{(1^m)}=\frac{m}2-\eta,
\quad
u_{(m)}=-\,\frac{m}2-\eta.
$$
Observe that the function $F_\mu(u)$ has
no pole at $u=u_\mu$ when $\g=\soN$,$\mu=(m)$ or $\g=\spN$,$\mu=(1^m)$.
Further, when $\g=\soN$,$\mu=(1^m)$ or $\g=\spN$,$\mu=(m)$
the order of the pole of $F_\mu(u)$ at $u=u_\mu$ does not exceed~1.
So the following definition 
of the normalizing factor $\phi_\mu(u)\in\CC(u)$
makes the product $\phi_\mu(u)\cdot F_\mu(u)$ regular at $u=u_\mu\tss$:
$$
\phi_{(1^m)}(u)
=
\cases
(u-\tsize\frac{m}2+\frac12)/(u+\frac12)&\quad\text{for}\ \ \g=\soN\,,
\\
1&\quad\text{for}\ \ \g=\spN\,;
\endcases
\Tag{3.333333}
$$

\kern-20pt

$$
\,\,\phi_{(m)}(u)
=
\cases
1&\quad\text{for}\ \ \g=\soN\,,
\\
(u+\tsize\frac{m}2-\frac12)/(u-\frac12)&\quad\text{for}\ \ \g=\spN\,.
\endcases
\Tag{3.333333333}
$$
Let $\tr$ be the standard trace on $\EndCN^{\ot m}$.
The following twin theorems constitute the main result of this section\tss;
\text{cf. [\tss O1\tss,\tss Theorem 1.3\tss]
and [\tss N2\tss,\tss Theorem 5.3\tss].}

\vbox{
\proclaim{Theorem 3.2}
\!For $\mu=(1^{2k})$ \!the value of
$\phi_\mu(u)\cdot(\tr\ot\id)\bigl(F_\mu(u)\bigr)$ 
\text{\!at $u=u_\mu$ \!is $\!C_k\!$.}
\endproclaim

\proclaim{Theorem 3.3}
\!For $\mu=({2k})$ \!the value of
$\phi_\mu(u)\cdot(\tr\ot\id)\bigl(F_\mu(u)\bigr)$ 
\text{\!at $u=u_\mu$ \!is $\!D_k\!$.}
\endproclaim
}

\nt
We will now present the proofs of both theorems as a chain
of propositions.

\proclaim{Proposition 3.4}
For each $\mu=(1^m)\tss,(m)$ we have
$(\tr\ot\id)\bigl(F_\mu(u)\bigr)\in\,\Zg(u)$.
\endproclaim

\demo{Proof}
Regard the group $G$ as a subgroup in $GL_N\subset\EndCN$.
Consider the adjoint action $\ad$ of the group $G$ in the enveloping algebra
$\Ug$.
Observe that by the definition \(3.00) for any element $g\in G$ we have
$$
(\id\ot\ad g)\,(F)=(g\ot1)\cdot F\cdot(g\1\ot1)\,.
$$
Each of the elements $A_m$,$B_m$ and $Q_{pq}$
in $\EndCN^{\ot m}$ commutes with $g^{\ot m}$. Thus
$$
(\id\ot\ad g)\,\bigl(F_\mu(u)\bigr)=
\bigl(g^{\ot m}\ot1\bigr)\cdot F_\mu(u)\cdot\bigl((g\1)^{\ot m}\ot1\bigr)\,.
$$
Hence
$$
(\tr\ot\ad g)\,\bigl(F_\mu(u)\bigr)=(\tr\ot\id)\,\bigl(F_\mu(u)\bigr)
\quad\square
$$
\enddemo

\nt
The proof of the next proposition is based on the following
simple lemma. Consider the product in \text{$\EndCN^{\ot3}(u,v,w)$}
$$
R_{12}(u,v)\,\Rt_{13}(u,w)\,\Rt_{23}(v,w)=
\Rt_{23}(v,w)\,\Rt_{13}(u,w)\,R_{12}(u,v)\,.
\Tag{3.03}
$$
The equality in \(3.03) can be verified directly. Notice
that the factor $\Rt_{23}(v,w)$ at either side of this equality
has a pole at $v+w=0$.
Yet here is the simple lemma.

\proclaim{Lemma 3.5}
The restrictions of \(3.03) to the two sets of $(u,v,w)$
where $v=u\pm1$, are regular at $v+w=0$.
\endproclaim

\demo{Proof}
Due to the relation $Q_{13}\,Q_{23}=P_{12}\,Q_{23}$ 
we have the equality
$$
R_{12}(u,u\pm1)\,\Rt_{13}(u,w)\,\Rt_{23}(u\pm1,w)=
\bigl(1\pm P_{12}\bigr)\cdot\biggl(1+\frac{Q_{13}+Q_{23}}{u+w}\biggr)\,.
$$
The latter rational function of $u,w$ is manifestly regular at $w=-u\mp1$
\enddemos 

\nt 
Let us identify the elements $A_{m-1}$ and $B_{m-1}$
of the algebra 
$\EndCN^{\ot\tss(m-1)}$ with their images under the natural embedding 
$$
\EndCN^{\ot\tss(m-1)}\to\EndCN^{\ot m}:\tss X\mapsto X\ot1\tss.
$$
In the next proposition we as always arrange non-commuting
factors corresponding to the indices $p=1\lc m-1$ from left to right.

\proclaim{Proposition 3.6}
We have the equalities in $\EndCN^{\ot m}(u)$
$$
\align
A_{m-1}\cdot\prod_{p=1}^{m-1}\,\Rt_{pm}(u-p+1,u-m+1)
&=A_{m-1}\cdot\biggl(1+\frac{Q_{1m}+\ldots+Q_{m-1,m}}{2u-m+1}\biggr)\tss,
\\
B_{m-1}\cdot\prod_{p=1}^{m-1}\,\Rt_{pm}(u+p-1,u+m-1)
&=B_{m-1}\cdot\biggl(1+\frac{Q_{1m}+\ldots+Q_{m-1,m}}{2u+m-1}\biggr)\tss.
\endalign
$$
\endproclaim

\demo{Proof}
We will give the proof of the first equality, the proof
of the second is quite similar.
Denote by $A(2u)$ the rational function of $2u$ at the left hand side
of the first relation. The value of this function at $2u=\infty$ is $A_{m-1}$.
Moreover, its residue at $2u=m-1$ is $Q_{1m}+\ldots+Q_{m-1,m}$.
So to prove the first relation it suffices to show that $A(2u)$ has no poles
except for a simple pole at $2u=m-1$.
Let an index $p\in\{\tss2\lc m-1\tss\}$ be fixed.
The corresponding factor $\Rt_{pm}(u-p+1,u-m+1)$ in $A(2u)$
has a pole at $2u=p+m-2$.
But when estimating from above the order of the pole at $2u=p+m-2$ of $A(2u)$,
that factor does not count. Indeed, the element $A_{m-1}\in\EndCN^{\ot m}$
is divisible on the right by $R_{p-1,p}(u-p+2,u-p+1)$.
Any factor in $A(2u)$ corresponding to an index less than $p-1$ commutes
with $R_{p-1,p}(u-p+2,u-p+1)$. But the product
$$
R_{p-1,p}(u\!-\!p+\!2,u\!-\!p+\!1)\,
\Rt_{p-1,m}(u\!-\!p+\!2,u\!-\!m\!+\!1)\,
\Rt_{pm}(u\!-\!p+\!1,u\!-\!m\!+\!1)
$$
has no pole at $2u=p+m-2$ due to Lemma 3.5
\enddemos

\nt
We will make use of the
well known decompositions of the elements \text{$A_m$ and $B_m$}
$$
\align
A_{m}\tss
&=\frac1{m\tss!}\tss\cdot\tss\prod_{p=1}^{m-1}\,
R_{pm}(1-p,1-m)\tss\ldots\tss R_{p,p+1}(1-p,-\,p)\tss,
\Tag{3.04}
\\
B_{m}\tss
&=\frac1{m\tss!}\tss\cdot\tss\prod_{p=1}^{m-1}\,
R_{pm}(p-1,m-1)\tss\ldots\tss R_{p,p+1}(p-1,+\,p)
\Tag{3.05}
\endalign
$$
in $\EndCN^{\ot m}$. They can be easily verified by the induction~on~$m$;
see for instance [\tss MNO\tss,\tss Section 2.3\tss].
Next proposition provides an alternative definition for each of
the functions $F_{(1^m)}(u)$ and $F_{(m)}(u)$.
\text{This definition was motivated by [\tss C1\tss,\tss S2\tss].}

\proclaim{Proposition 3.7}
We have the equalities in $\EndCN^{\ot m}\ot\Ug(u)$
$$
\align
F_{(1^m)}(u)
&=(A_m\ot1)\,\cdot\,\prod_{q=1}^m\
\Bigl(\ \prod_{p=1}^{q-1}\ \Rt_{pq}(u-p+1,u-q+1)\Bigr)
\,F_q(u-q+1)
\tss,
\\
F_{(m)}(u)
&=(B_m\ot1)\,\cdot\,\prod_{q=1}^m\
\Bigl(\ 
\prod_{p=1}^{q-1}\ \Rt_{pq}(u+p-1,u+q-1)
\Bigr)
\,F_q(u+q-1)\tss.
\endalign
$$
\endproclaim

\demo{Proof}
It suffices to show that at the right hand sides of the above formulas
for the rational functions
$F_{(1^m)}(u)$ and $F_{(m)}(u)$ one can replace the ordered products
$$
\prod_{p=1}^{m-1}\,\Rt_{pm}(u-p+1,u-m+1)
\quad\text{and}\quad
\prod_{p=1}^{m-1}\,\Rt_{pm}(u+p-1,u+m-1)
$$
by 
$$
1+\frac{Q_{1m}+\ldots+Q_{m-1,m}}{2u-m+1}
\quad\ \ \text{and}\ \ \quad
1+\frac{Q_{1m}+\ldots+Q_{m-1,m}}{2u+m-1}
$$
respectively. Using the decompositions \(3.04),\(3.05)
for $m-1$ instead of $m$ 
and the relations \(3.111),\(3.03) one can insert 
respectively $A_{m-1}$ and $B_{m-1}$ right after
the factors $F_{m-1}(u-m+2)$ and $F_{m-1}(u+m-2)$ in the above
formulas for $F_{(1^m)}(u)$ and $F_{(m)}(u)$. 
Here we also use the equalities
$A_{m-1}^{\tss2}=A_{m-1}$ and $B_{m-1}^{\tss2}=B_{m-1}$.
Now the required statement follows directly from Proposition 3.6
\enddemos

\nt
The main part of our proof of Theorems 3.2 and 3.3 will employ
Corollary 2.5. As in \text{Section~2,} let $\la$ be any partition with
$\ell(\la)\le n$. Let $V_\la$ be the corresponding irreducible
$\g$-module. Put $l=|\la|$. We will use the classical realization 
of $V_\la$ as a submodule in $(\CC^N)^{\ot\tss l}$ consisting of
traceless tensors. Recall that a vector
$\xi\in(\CC^N)^{\ot\tss l}$ is called
{\it traceless} if $Q_{rs}\cdot\xi=0$ whenever $1\le r<s\le l$. The
subspace in $(\CC^N)^{\ot\tss l}$ formed by all such vectors will
be denoted by $V$.

Let $U_\la$ be the irreducible module over the Lie algebra $\glN$
corresponding to the partition $\la$. Fix any embedding of the $\glN$-module 
$U_\la$ to $(\CC^N)^{\ot\tss l}$. The subspace $U_\la\cap V$
in $U_\la$ is preserved by 
the action of the subalgebra $\g\subset\glN$. For
$\g=\so_{2n+1}\tss,\tss\sp_{2n}$ this subspace is isomorphic to $V_\la$ as
$\g$-module. For $\g=\so_{2n}$ it is isomorphic to $V_\la$ only
if $\la_n=0$. Otherwise $U_\la\cap V$ splits into direct sum of the
$\so_{2n}$-modules $V_\la$ and $V_\las$.
All these statements are contained in \text{[\tss W,\tss Section V.C\tss].}

Our proof of Theorems 3.2 and 3.3 will be based on another simple lemma.
Put
$$
E(u)\,=-\,u+\sum_{ij}\,E_{ij}\ot E_{ji}\,\in\,
\EndCN\ot\UglN\tss[\tss u\tss]\tss.
$$
Let
$$
\Et(u)\,=-\,u+\sum_{ij}\,\ep_{ij}\cdot E_{ij}\ot E_{-i,-j}
\,\in\,\EndCN\ot\UglN\tss[\tss u\tss]
$$
be the element
obtained from $E(u)$ by \text{applying} the transposition with respect
to the bilinear form $\langle\,\,,\,\rangle$ in the tensor factor $\EndCN$.
Then consider the element
$$
\qquad\qquad
\frac{\Et(\eta-u)\,E(\eta+u)}{u-\eta}\,\in\,\EndCN\ot\UglN\tss(\tss u\tss)
\tss.
\Tag{3.000222}
$$
We have the standard representation $\UglN\to\EndCN^{\ot\tss l}$ so 
the element \(3.000222) acts in the space $\CC^N\ot(\CC^N)^{\ot\tss l}$.
The element $F(u)\in\EndCN\ot\Ug\tss[\tss u\tss]$ also acts in the
space $\CC^N\ot(\CC^N)^{\ot\tss l}$ and the latter action preserves the
subspace $\CC^N\ot V$. 
Here is the lemma; cf. [O] and [\tss MNO\tss,\tss Section 3.5\tss].

\proclaim{Lemma 3.8}
Action of the element $F(u)$ in $\CC^N\ot V\!$ coincides with
that of \,\(3.000222).
\endproclaim

\demo{Proof}
Consider the image of $F(u)\in\EndCN\ot\Ug\tss[\tss u\tss]$
in $\EndCN^{\ot\tss(l+1)}\tss[\tss u\tss]$
under the
representation $\Ug\to\EndCN^{\ot\tss l}$.
By \(3.00) 
and \(3.025) this image equals
$$
P_{12}+\ldots+P_{1,l+1}-Q_{12}-\ldots-Q_{1,l+1}-u-\eta
\tss.
\Tag{3.06}
$$
On the other hand, the image of the element \(3.000222) in
$\EndCN^{\ot\tss(l+1)}\tss(\tss u\tss)$ under the
representation $\UglN\to\EndCN^{\ot\tss l}$ is
$$
\biggl(1+\frac{Q_{12}+\ldots+Q_{1,l+1}}{u-\eta}\biggr)
\cdot
\bigl(P_{12}+\ldots+P_{1,l+1}-u-\eta\tss\bigr)\tss.
\Tag{3.07}
$$
Using the relation \(3.01) and definition of $\eta$
one shows that the product \(3.07) equals \(3.06) plus the sum
$$
\sum_{r\neq s}\ \,\frac{Q_{1,r+1}\,P_{1,s+1}}{u-\eta}\,.
$$
Since $Q_{1,r+1}\,P_{1,s+1}=P_{1,s+1}\,Q_{s+1,r+1}$ for $r\neq s$,
the restriction of the latter sum to $\CC^N\ot V$ vanishes identically
\enddemos

\nt
Similarly to \(3.0003), for $p=1\lc m$ denote respectively by
$E_p(u)$,$\Et_p(u)$ the images in $\EndCN^{\ot m}\ot\UglN\,[\tss u\tss]$
of $E(u)$,$\Et(u)$ under the embedding $\io_p\ot\id$.

\proclaim{Proposition 3.9}
We have the relations in the algebra $\EndCN^{\ot2}\ot\UglN\,(u,v)$
$$
\align
R(u,v)\,E_1(u)\,E_2(v)&=E_2(v)\,E_1(u)\,R(u,v)\,,
\Tag{3.81}
\\
R(u,v)\,\Et_1(-\,u)\,\Et_2(-\,v)&=\Et_2(-\,v)\,\Et_1(-\,u)\,R(u,v)\,,
\Tag{3.82}
\\
\Et_1(-u)\,\Rt(u,v)\,E_2(v)&=E_2(v)\,\Rt(u,v)\,\Et_1(-\,u)\,.
\Tag{3.83}
\\
\hskip5cm&
\endalign
$$
\endproclaim

\kern-0.75cm
\demo{Proof}
As well as \(3.111), the relation \(3.81) can be verified directly.
Relation \(3.82) is obtained from \(3.81) by applying in the tensor
factor $\UglN$ the automorphism $E_{ij}\mapsto-\tss\ep_{ij}\cdot E_{-j,-i}$.
By applying to \(3.81) the transposition with respect to
$\langle\,\,,\,\rangle$ in the fisrt tensor factor of $\EndCN^{\ot2}$
we get \(3.83)
\enddemos

\nt
We will need two more propositions.
The first of them implies the classical Capelli
identity [\tss C\tss]. 
The second proposition implies the identity
from [\tss N1\tss,\tss Example 2\tss] mentioned in Section~1.
Observe that the $m$-th Capelli element of $\ZglN$ is the image
under the map $\tr\ot\id$ of the element from $\EndCN^{\ot m}\ot\UglN$
$$
(A_m\ot1)\cdot E_1(0)\,E_2(-1)\,\ldots\,E_m(1-m)\tss,
\nopagebreak
\Tag{3.000444}
$$
see \(0.3)\tss. Similarly, the element \(0.31) with $k=m$ is the image
under $\tr\ot\id$ of
$$
(B_m\ot1)\cdot E_1(0)\,E_2(1)\,\ldots\,E_m(m-1)\tss.
\Tag{3.0004445}
$$
In proofs of both propositions \text{we follow [M2].}
For further details see \text{[\tss N2\tss,\tss Section 5\tss].}

\proclaim{Proposition 3.10}
The element \(3.000444) of $\EndCN^{\ot m}\ot\UglN$
vanishes in every representation $\UglN\to\End(U_\la)$ where $\ell(\la)<m$.
\endproclaim

\proclaim{Proposition 3.11}
The element \(3.0004445) of $\EndCN^{\ot m}\ot\UglN$
vanishes in every representation $\UglN\to\End(U_\la)$ where $\la_1<m$.
\endproclaim

\demo{Proofs}
We give only the proof of Proposition 3.10.
The proof of Proposition~3.11 is similar and will be omitted.
We will show by induction on $m$ that the image in
$\EndCN^{\tss(m+l)}$ 
of \(3.000444) under the representation
$\UglN\to\EndCN^{\ot\tss l}$ equals
$$
A_m\,\cdot\!\sum_{r_1\ldots r_m}P_{1,m+r_1}\ldots P_{m,m+r_m}
\Tag{3.000555}
$$
where all the indices $r_1\lc r_m\in\{\tss 1\lc l\tss\}$
are pairwise distinct. The sum in \(3.000555) vanishes when
$l<m$. When $l\ge m$ but $\ell(\la)<m$ the required vanishing property 
will follow from the Young decomposition of the $\glN$-module
$U_{(1^m)}\ot(\CC^N)^{\ot\tss(l-m)}$
into irreducible components. 

By definition the image in $\EndCN^{\tss(m+l)}[\tss u\tss]$
of $E_p(u)\in\EndCN^{\ot m}\ot\UglN$ under
the representation $\UglN\to\EndCN^{\ot\tss l}$ is
$$
-u+\sum_r\, P_{p,m+r}
$$
where the sum is taken over $r=1\lc l$.
In particular, this provides the base $m=1$ for our induction.
Now suppose that $m>1$.
By the inductive assumption
the image in $\EndCN^{\ot\tss(m+l)}$ of the element \(3.000444) equals
$$
\gather
A_m\,\cdot\!\!\!\!\sum_{r_1\ldots r_{m-1}}\!\!
P_{1,m+r_1}\ldots P_{m-1,m+r_{m-1}}\,\cdot\,
\sum_r\, P_{m,m+r}\,\,+
\Tag{3.000666}
\\
(m-1)\,\cdot\, 
A_m\,\cdot\!\!\!\!\sum_{r_1\ldots r_{m-1}}\!\!
P_{1,m+r_1}\ldots P_{m-1,m+r_{m-1}}
\endgather
$$
where all the indices $r_1\lc r_{m-1}\in\{\tss 1\lc l\tss\}$
are pairwise distinct. Here the expression in the second line 
can be replaced by 
$$
-\,\,A_m\,\cdot\,\sum_{p=1}^{m-1}\,P_{pm}\,
\cdot\!\!\!\!\sum_{r_1\ldots r_{m-1}}\!\!
P_{1,m+r_1}\ldots P_{m-1,m+r_{m-1}}\,
$$
But the latter expression can be rewritten as 
$$
-\,\,A_m\,\cdot\!\!\!\!\sum_{r_1\ldots r_{m-1}}\!\!
P_{1,m+r_1}\ldots P_{m-1,m+r_{m-1}}\,\cdot\,
\sum_{p=1}^{m-1}\,P_{m,m+r_p}\,.
$$
Hence after this replacement the sum \(3.000666) becomes
evidently equal to \(3.000555)
\enddemos
\nopagebreak
\vbox{
\nt
Note that Proposition 3.10 can be reformulated in the following way.

\proclaim{Corollary 3.12}
The element of $\EndCN^{\ot m}\ot\UglN$
$$
(A_m\ot1)\cdot\Et_1(1-m)\,\ldots\,\Et_{m-1}(-1)\,\Et_m(0)
\nopagebreak
\Tag{3.0004441}
$$
vanishes in every representation $\UglN\to\End(U_\la)$ where $\ell(\la)<m$.
\endproclaim
}

\demo{Proof}
By using the decomposition \(3.04) and the relation \(3.81) the element
\(3.000444) of the algebra $\EndCN^{\ot m}\ot\UglN$ can be rewritten as
$$
E_m(1-m)\,\ldots\,E_2(-1)\,E_1(0)\cdot(A_m\ot1)\tss.
\Tag{3.0004442}
$$
Applying to \(3.0004442)
the transposition in each tensor factor of $\EndCN^{\ot m}$
relative to $\langle\,\,,\,\rangle$ we get
$$
(A_m\ot1)\cdot\Et_m(1-m)\,\ldots\,\Et_2(-1)\,\Et_1(0)\tss.
$$
Conjugate the latter element by $P_m\ot 1,$ where $P_m\in\EndCN^{\ot m}$
corresponds to the operator
$$
e_{i_1}\ot e_{i_2}\ot\ldots\ot e_{i_m}\mapsto
e_{i_m}\ot\ldots\ot e_{i_2}\ot e_{i_1}
$$
in $(\CC^N)^{\ot m}$. Then we obtain exactly the element \(3.0004442)
\enddemos 

\nt
By using the decomposition \(3.05) and again Proposition 3.9,
we obtain a corollary to Proposition 3.11 parallel to Corollary 3.12.

\proclaim{Corollary 3.13}
The element of $\EndCN^{\ot m}\ot\UglN$
$$
(B_m\ot1)\cdot\Et_1(m-1)\,\ldots\,\Et_{m-1}(1)\,\Et_m(0)
$$
vanishes in every representation $\UglN\to\End(U_\la)$ where $\la_1<m$.
\endproclaim

\nt
If we replace each appearance of the function $F(u)$ 
in the definitions of $F_{(1^m)}(u)$ and $F_{(m)}(u)$ by 
\(3.000222), we obtain the functions valued in
$\EndCN^{\ot\tss m}\ot\UglN$
$$
\gather
(A_m\ot1)\cdot\prod_{q=1}^m\,
\biggl(\! 
\Bigl(\tss1+\frac{Q_{1q}+\ldots+Q_{q-1,q}}{2u-q+1}\tss\Bigr)\ot1  
\tss\!\biggr)
\,
\frac{\Et_q(\eta\!-\!u\!+\!q\!-\!1)\tss
 E_q(\eta\!+\!u\!-\!q\!+\!1)}{u-q+1-\eta}\tss,
\\
(B_m\ot1)\cdot\prod_{q=1}^m\,
\biggl(\!
\Bigr(\tss1+\frac{Q_{1q}+\ldots+Q_{q-1,q}}{2u+q-1}\tss\Bigr)\ot1
\tss\!\biggr)
\,
\frac{\Et_q(\eta\!-\!u\!-\!q\!+\!1)\tss
E_q(\eta\!+\!u\!+\!q\!-\!1)}{u+q-1-\eta}\tss.
\endgather
\nopagebreak
$$
We will denote these rational functions by $G_{(1^m)}(u)$ and
$G_{(m)}(u)$ respectively.

Now for either of the two partitions $\mu=(1^{2k}),(2k)$ of $m=2k$
\text{take the value of}
$$
\phi_\mu(u)\cdot F_\mu(u)\in\EndCN^{\ot m}\ot\Ug(u)
$$
at $u=u_\mu$.
Let us consider the image of this value with respect to the
representation $\Ug\to\End(V_\la)$.
For $\mu=(1^{2k})$ we will show that this image is zero
when $\ell(\la)<k$. For $\mu=({2k})$ we will show  that this image is zero
when $\la_1<k$. Thanks to Lemma 3.8 it
suffices to prove the next twin propositions.

\proclaim{Proposition 3.14}
For $\mu=(1^{2k})$ the value of the function $\phi_\mu(u)\cdot G_\mu(u)$
at $u=u_\mu$ vanishes in every representation $\UglN\to\End(U_\la)$
with $\ell(\la)<k$.
\endproclaim
\nopagebreak
\proclaim{Proposition 3.15}
For $\mu=({2k})$ the value of the function $\phi_\mu(u)\cdot G_\mu(u)$
at $u=u_\mu$ vanishes in every representation $\UglN\to\End(U_\la)$
with $\la_1<k$.
\endproclaim

\demo{Proofs}
We will give the proof of Proposition 3.14\tss; the proof of
\text{Proposition 3.15} is similar and will be omitted. We will prove that
the function $G_{(1^{2k})}(u)$ equals
$$
\gather
(A_{2k}\ot1)\,\cdot\,\prod_{q=1}^{2k}\,\tss\Et_q(\eta-u+q-1)\,\,\x 
\Tag{3.1313}
\\
\prod_{q=1}^{2k}\,
\biggl(\!
\Bigl(\tss1+\frac{Q_{1q}+\ldots+Q_{q-1,q}}{2u-q+1}\tss\Bigr)\ot1
\tss\!\biggr)
\,\cdot\,
\prod_{q=1}^{2k}\,\frac1{u-q+1-\eta}\,\,\x
\\
(A_{2k}\ot1)\,\cdot\,\prod_{q=1}^{2k}\,\tss E_q(\eta+u-q+1)\quad  
\endgather
\nopagebreak
$$
where the non-commuting
factors corresponding to the indices $q=1\lc 2k$ are as usual 
arranged from left to right. 
\text{\hskip-0.3pt Note that we have defined $u_{(1^{2k})}=k-\eta$.} 
But in the third line of 
\(3.1313) we get the function valued in $\EndCN^{\ot\tss 2k}\ot\UglN$
$$
(A_{2k}\ot1)\,\cdot\,
E_1(u+\eta)\ldots E_k(u+\eta-k+1)\,\cdot\,
E_{k+1}(u+\eta-k)\ldots E_{2k}\bigl(u+\eta-2k+1)\,.
$$
Its value at $u=k-\eta$ vanishes 
in every representation $\UglN\to\End(U_\la)$ with $\ell(\la)<k$.
To see this consider the last $k$ factors in the above product and
apply Proposition 3.10 to $m=k$.

Suppose that $\g=\soN$.
Then $\eta=\frac12$ and $u_{(1^{2k})}=k-\frac12$.
Observe that $Q_{pq}\cdot A_{2k}=0$ for any $1\le p<q\le 2k$.
Therefore the first product in the second line of \(3.1313) can
be replaced by $1$. The second product in that line has a pole
at $u=k-\frac12$ of degree $1$.
Thus we will obtain Proposition 3.14 for $\g=\soN$ since in this case
$$
\phi_{(1^{2k})}(u)=\tsize(u-k+\frac12)/(u+\frac12)\tss.
$$

Now suppose that $\g=\spN$.
Then $\eta=-\,\frac12$\,,\,$u_{(1^{2k})}=k+\frac12$ while
$\phi_{(1^{2k})}(u)=1$.
If $k=1$ then the rational function in second line of
\(3.1313) has no pole at \text{$u=\frac32$\tss;}
so we get Proposition 3.14
by considering the product in the third line of \(3.1313) only. 
Assume that $k\ge2$.
Then the function in the second line of
\(3.1313) has a pole at $u=k+\frac12$ of degree at most 1. 
But in the first line of
\(3.1313) we then get
$$
\align
(A_{2k}\ot1)\,\cdot\,
&\Et_1(-\,\tsize\frac12-u)\,\Et_2(\frac12-u)\cdot
\Et_3(\frac32-u)\,\ldots\,\Et_{k+2}(k+\frac12-u)\,\,\x
\Tag{3.131313}
\\
\qquad\qquad\qquad\qquad\quad
&\tss\qquad
\Et_{k+3}(k+\tsize\frac32-u)\,\ldots\,\Et_{2k}(2k-\frac32-u)\,.
\endalign
$$
The value of \(3.131313) at $u=k+\frac12$ vanishes 
in any representation \text{$\UglN\to\End(U_\la)$} with $\ell(\la)<k$.
To see this take the last $k$ factors in the first line of
\(3.131313) and apply 
Corollary 3.12 to $m=k$.
Thus we will obtain Proposition 3.14 for $\g=\spN$.

It remains to show that the function $G_{(1^{2k})}(u)$ indeed equals
\(3.1313). By using \(3.04) and \(3.82) we obtain for any $m$ the relation 
in $\EndCN^{\ot m}\ot\UglN(u)$
$$
(A_{m}\ot1)\,\cdot\,\Et_1(u)\,\ldots\,\Et_m(u+m-1)=
\Et_m(u+m-1)\,\ldots\,\Et_1(u)\,\cdot\,(A_{m}\ot1)\,.
$$
In particular, the expression
in the first line of \(3.1313) is divisible on the right by $A_{2k}\ot1$.
But due to Proposition 3.6 we have
$$
\gather
A_{2k}\cdot\prod_{q=1}^{2k}\,
\Bigl(\tss1+\frac{Q_{1q}+\ldots+Q_{q-1,q}}{2u-q+1}\tss\Bigr)=
A_{2k}\cdot
\prod_{q=1}^{2k}\,A_{q-1}\,
\Bigl(\tss1+\frac{Q_{1q}+\ldots+Q_{q-1,q}}{2u-q+1}\tss\Bigr)
\\
=\,A_{2k}\cdot
\prod_{q=1}^{2k}\
\Bigl(\ \prod_{p=1}^{q-1}\,\Rt_{pq}(u-p+1,u-q+1)\tss\Bigr)\,.
\endgather
$$
So \(3.1313) equals
 
\kern-25pt
$$
\align
(A_{2k}\ot1)\,&\cdot\,\prod_{q=1}^{2k}\,\tss\Et_q(\eta-u+q-1)\,\,\x 
\\
\qquad
\prod_{q=1}^{2k}\
\Bigl(\ 
\prod_{p=1}^{q-1}\,\Rt_{pq}(u-p&+1,u-q+1)\tss\Bigr)
\,\cdot\,
\prod_{q=1}^{2k}\,\frac1{u-q+1-\eta}\,\,\x
\\
(A_{2k}\ot1)\,&\cdot\,\prod_{q=1}^{2k}\,\tss E_q(\eta+u-q+1)\,.\quad  
\endalign
$$
By using only \(3.83) we can transform this expression 
to
$$
(A_{2k}\ot1)\cdot
\prod_{q=1}^{2k}\ 
\biggl(\ 
\prod_{p=1}^{q-1}\,\Rt_{pq}(u\!-\!p\!+\!1,u\!-\!q\!+\!1)
\biggr)
\,
\frac{\Et_q(\eta\!-\!u\!+\!q\!-\!1)\tss E_q(\eta\!+\!u\!-\!q\!+\!1)}
{u\!-\!q\!+\!1\!-\!\eta}\,
\,.
$$
The latter expression equals $G_{(1^{2k})}(u)$. To prove this observe that
the element \(3.000222) satisfies a relation in $\EndCN^{\ot2}\ot\UglN(u,v)$
similar to \(3.111)
$$
\align
&R(u,v)\cdot\Et_1(\eta-u)\,E_1(\eta+u)\cdot
\Rt(u,v)\cdot\Et_2(\eta-v)\,E_2(\eta+v)=
\\
&\Et_2(\eta-v)\,E_2(\eta+v)\cdot\Rt(u,v)\cdot
\Et_1(\eta-u)\,E_1(\eta+u)\cdot R(u,v)\,.
\endalign
\nopagebreak
$$
This relation can be easily deduced from \(3.81),\(3.82),\(3.83).
So we can repeat the arguments from Proposition 3.7
for the function \(3.000222) instead of $F(u)$
\enddemos

\nt
Now for $\mu=(1^{2k})$ take the value of
the function $\phi_\mu(u)\cdot(\tr\ot\id)\bigl(F_\mu(u)\bigr)$
at $u=u_\mu$. Denote this value by $Z$. By Proposition 3.1 we have
$Z\in\Zg$. The degree of $Z$ with respect to the standard
filtration on $\Ug$ does not exceed $2k$ by the definition of $F_\mu(u)$.
By Corollary 2.5 we obtain from Proposition 3.14 that $Z$ equals $C_k$
up to a scalar multiple.
To complete the proof of Theorem 3.2 it
remains to show that this multiple is actually $1$.
Similarly, for $\mu=(2k)$ due to Proposition 3.15 the value of
$\phi_\mu(u)\cdot(\tr\ot\id)\bigl(F_\mu(u)\bigr)$ at $u=u_\mu$
is $D_k$ up to a scalar multiple. So to complete the proof of Theorem 3.3
\text{it remains to show that this multiple is also $1$.}

For $\mu=(1^{2k}),(2k)$ take the image
under the Harish-Chandra isomorphism
$$
\om\tss\Bigl(\,\phi_\mu(u)\cdot(\tr\ot\id)\bigl(F_\mu(u)\bigr)\Bigr)\tss.
\Tag{3.333}
$$
Consider this image as a polynomial in $\la_1\lc\la_n$ with
the coefficients from $\CC(u)$. The next two propositions
complete the proofs of \text{Theorems 3.2 \!and \!3.3 respectively.}
They also justify our choice \(3.333333),\(3.333333333)
of the normalizing factor $\phi_\mu(u)\in\CC(u)$.

\proclaim{Proposition 3.16}
For $\mu=(1^{2k})$ the leading component of the polynomial \(3.333) 
is $(-1)^k\tss e_k(\la_1^2\lc\la_n^2)$.
\endproclaim

\proclaim{Proposition 3.17}
For $\mu=(2k)$ the leading component of the polynomial \(3.333)
is $h_k(\la_1^2\lc\la_n^2)$.
\endproclaim

\demo{Proofs}
We shall give the proof of Proposition 3.16\tss;
the proof of Proposition~3.17 is similar and will be omitted.
We will write $m=2k$ for short.
First let us prove that leading
component of the polynomial \(3.333) with $\mu=(1^m)$
does not depend on the parameter $u$.
We will use the equality in 
\text{$\EndCN^{\ot m}\ot\Ug(u)$}
$$
F_{(1^m)}(u)=(A_m\ot1)\,\tss\cdot\,\prod_{q=1}^m\ F_q(u-q+1)\tss
\biggl(\!\hskip-1pt
\Bigl(\tss1+\frac{Q_{q,q+1}+\ldots+Q_{qm}}{2u-2q+1}\tss\Bigr)\ot1
\tss\!\biggr)
$$
where the factors corresponding to the indices $q=1\lc m$ are as usual
arranged from left to right. It can be proved by using
Proposition 3.7, see also Proposition~3.6.
Denote by $A_{m-1}^{\tss\prime}$ the image of the element $A_{m-1}$
with respect to the embedding
$$
\EndCN^{\ot\tss(m-1)}\to\EndCN^{\ot\tss m}:\tss\tss
X\mapsto 1\ot X\tss.
\Tag{3.*0}
$$
Let us use the well known decomposition 
$A_m=T_m\cdot A_{m-1}^{\tss\prime}$
in $\EndCN^{\ot\tss m}$ where
$$
T_m=
\bigl(\tss1-P_{12}-\ldots-P_{1m}\bigr)/m\tss.
$$
By this decomposition and by the above formula for $F_{(1^m)}(u)$
that function equals
$$
\bigl(T_m\ot1\bigr)
\,F_1(u)\,
\biggl(\!
\Bigl(\tss1+\frac{Q_{12}+\ldots+Q_{1m}}{2u-1}\tss\Bigr)\ot1\!\tss
\biggr)
\,F_{(1^{m-1})}^{\tss\prime}(u-1)
\Tag{3.*1}
$$
where $F_{(1^{m-1})}^{\tss\prime}(u)$ stands
for the image of 
$F_{(1^{m-1})}(u)\in\EndCN^{\ot\tss(m-1)}\ot\Ug(u)$
with respect to the embedding \(3.*0).
Take the standard filtration on the algebra $\Ug$.
Extend it to the algebra $\Ug(u)$ by setting $\deg u=0$.
Denote
$$
\psi\tss(u)=
\cases
(u+\tsize\frac12)/(u-\frac12)&\quad\text{for}\ \ \g=\soN\,,
\\
1&\quad\text{for}\ \ \g=\spN\,.
\endcases
$$
Then the image in $\Ug(u)$ of the element \(3.*1) under $\tr\ot\id$
is equal to that of
$$
\psi\tss(u)\cdot
\bigl(T_m\ot1\bigr)
\,F_1(u)\,
F_{(1^{m-1})}^{\tss\prime}(u-1)
\Tag{3.*3}
$$
plus lower degree terms.
Indeed, due to \(3.00) and \(3.025) for any $r=2\lc m$ the element 
$$
\bigl(P_{1r}\ot1\bigr)\cdot F_1\cdot\bigl(Q_{1r}\ot1\bigr)
\in\EndCN^{\ot m}\ot\Ug
$$
is equal to $-\tss F_1\cdot\bigl(Q_{1r}\ot1\bigr)$ if $\g=\soN$ and to
$F_1\cdot\bigl(Q_{1r}\ot1\bigr)$ if $\g=\spN$.
Therefore
$$
\bigl(T_m\ot1\bigr)
\,F_1\cdot
\biggl(\!
\Bigl(\tss1+\frac{Q_{12}+\ldots+Q_{1m}}{2u-1}\tss\Bigr)\ot1\!\tss
\biggr)=
\psi\tss(u)\cdot
\bigl(T_m\ot1\bigr)
\,F_1
\nopagebreak
$$
plus a certain element of $\EndCN^{\ot m}\ot\Ug(u)$ which is
antisymmetric under the transposition in the first tensor
factor $\EndCN$
with respect to the bilinear form $\langle\,\,,\,\rangle$.
So the image of the latter element under the map $\tr\ot\id$ is zero.

Let us now use the definition of $F_{(1^m)}(u)$ directly. 
By the properties of trace, the image of
$F_{(1^m)}(u)$ under $\tr\ot\id$ will remain unchanged when we move $A_m\ot1$
to the rightmost position and then apply the conjugation by
$\bigl(P_{12}\ldots P_{1m}\bigr)\ot1$.
After that we will decompose $A_m=T_m\cdot A_{m-1}^{\tss\prime}$
and move $A_{m-1}^{\tss\prime}$ to the leftmost position.
The transformed expression for the element $F_{(1^m)}(u-1)$
\text{can be written as} 
$$
F_{(1^{m-1})}^{\tss\prime}(u-1)
\biggl(\!
\Bigl(\tss1+\frac{Q_{12}+\ldots+Q_{1m}}{2u-m-1}\tss\Bigr)\ot1
\biggr)
\,F_1(u-m)\,\bigl(T_m\ot1\bigr)\tss.
$$
But the image of the latter element under the map $\tr\ot\id$
is equal to the image of
$$
\psi\tss(u-\tsize\frac{m}2\dsize)\cdot
F_{(1^{m-1})}^{\tss\prime}(u-1)\,
F_1(u-m)\,\bigl(T_m\ot1\bigr)
\Tag{3.*4}
$$
plus lower degree terms in $\Ug(u)$.
Note that by the definition \(3.333333)
$$
\phi_{(1^m)}(u)\cdot\psi\tss(u)=
\phi_{(1^m)}(u-1)\cdot\psi\tss(u-\tsize\frac{m}2)\tss.
$$
Therefore by comparing the expressions \(3.*3) and \(3.*4)
we obtain that \text{the degree} $m$ component of the polynomial \(3.333)
with $\mu=(1^m)$ does not change if we replace $u$ by $u-1$.
As this component is a rational
function in $u$,\, \text{it must be constant in $u$.}

To complete the proof of Proposition 3.16 it now suffices to
determine the leading component of 
\(3.333) with $\mu=(1^m)$ at $u\to\infty$.
This can be done as follows.
In the definition of $F_{(1^{m})}(u)$
replace each factor $F_q(u-q+1)$ by $F_q$ and note that
$\bigl(Q_{1q}+\ldots+Q_{q-1,q}\bigr)/(2u-q+1)\to0$
when  $u\to\infty$.
Since $\phi_{(1^{m})}(u)\to1$ when $u\to\infty$,
it now remains to calculate \text{the leading component of} the polynomial
$$
\om\tss\Bigl(
(\tr\ot\id)\tss\bigl(\tss(A_{m}\ot 1)\cdot F_1\ldots F_{m}\tss\bigr)
\!\tss\Bigr)
$$
with $m=2k$. This calculation is immediate, and yields
$(-1)^k e_k(\la_1^2\lc\la_n^2)$  
\enddemos


\kern6pt
\line{\bf4.\ Capelli identity for the Lie algebra $\soN$\hfill}
\section{\,}\kern-20pt

\vbox{\nt
In this section we will keep the positive integer $N$ fixed.
Let us now fix one more positive integer $m$.
We will consider a natural representation of the enveloping algebra
$\UsoN$ by differential operators with polynomial coefficients on the
space $\CMN$. We will give an explicit formula for the differential
operator corresponding to the element $C_k\in\ZsoN$
with \text{$k=1\lc n$.} We will show that
this formula may be regarded as an analogue for $\soN$
of the Capelli identity [\tss C\tss] for $\glN$.

As well as in Sections 2 and 3 here the indices $i$ and $j$ will run 
through the set $\{\tss-\tss n\lc\!-\!1,1\lc n\tss\}$ if $N=2n$ and 
through the set 
$\{\tss-\tss n\lc\!-\!1,0,1\lc n\tss\}$ if $N=2n+1$. The indices $a$ and $b$
will run through $\{\tss1\lc m\tss\}$. Let $x_{ai}$ be
independent commuting variables and $\d_{ai}$ be the corresponding
differentiations.
Denote by $\P$ the ring of polynomials in the variables
$x_{ai}$. Take the action \(4.0) of the Lie algebra $\glN$ on $\P$.
We will use the embedding of the group $O_N$ into $GL_N$ from
Section 2. Then restriction of the above action of $\glN$ in $\P$
\text{to $\soN$ is}
$$
F_{ij}\,\mapsto\,\sum_a\ \bigl(\,x_{ai}\d_{aj}-x_{a,-j}\d_{a,-i}\tss\bigr)\,.
\Tag{4.1}
$$  
}
Denote by $\PD$ the ring of differential operators in $x_{ai}$
with polynomial coefficients. Then \(4.1) determines a representation
$\ga:\UsoN\to\PD$. We will give a formula for the image in $\PD$ 
of the element $C_k\in\ZsoN$ with respect~to~$\ga$.

By definition we have $F_{i,-j}=-\tss F_{j,-i}$.
Let $I$ be any set of indices $i$ with $2k$ \text{elements.}
Suppose these elements are $i_1<\ldots<i_{2k}$.
\text{Take an element of $\UsoN$} 
$$
\Phi_I=\sum_{\si}\ \frac{\sgn(\si)}{2^k\tss k\tss!}\cdot
F_{i_{\si(1)},-i_{\si(2)}}\ldots F_{i_{\si(2k-1)},-i_{\si(2k)}}
\Tag{4.2}
$$
where $\si$ runs through the symmetric group $S_{2k}$.
The element \(4.2) is the {\it Pfaffian} of the
antisymmetric matrix 
$[\tss F_{i_p,-i_q}\tss]$ with the indices $p,q=1\lc 2k$.
Note that the factors of the monomial in the sum \(4.2) do not commute
in general. Still we will write
$$
\Phi_I=\Pf\,[\tss F_{i_p,-i_q}\tss]\,;
\ \quad p,q=1\lc 2k.  
$$

Let us denote $\Is=\{-\tss i_{2k}\lc\!-\!i_1\}$. 
Further, let $A$ be any set of indices $a$ with $k$ elements.
Suppose these elements are $a_1\lc a_k$. Then denote in $\PD$
$$
\Om_{AI}=\sum_J\,\sgn(JJ^\prime)\cdot
\det\tss[\tss x_{a_pj_q}\tss]\tss
\det\tss[\tss\d_{a_p,-j_q^{\tss\prime}}\tss]
\Tag{4.4}
$$
where the sum is taken over all partitions of $I$ into two subsets
$J=\{\tss j_1\lc j_k\tss \}$ and
$J^\prime=\{\tss j_1^{\tss\prime}\lc j_k^{\tss\prime}\tss\}$.
Here the determinants are taken with respect to the indices $p,q=1\lc k$
and $\sgn(JJ^\prime)$ is the sign of the permutation
$\bigl(\tss j_1,j_1^{\tss\prime}\lc j_k,j_k^{\tss\prime}\tss\bigr)$
of the sequence $(i_1\lc i_{2k})$.
Next theorem \text{provides an explicit formula for $\ga(C_k)$.}

\proclaim{Theorem 4.1}
We have the equalities respectively in $\UsoN$ and $\PD$
$$
C_k=(-1)^k\cdot\sum_I\,\Phi_I\tss\Phi_\Is
\ \quad\text{and}\ \quad
\ga(\Phi_I)=\sum_A\,\Om_{AI}\,.
\Tag{4.3}
$$
\endproclaim
 
\demo{Proof}
We start with verifying the second equality in \(4.3).
\text{By definition $\Phi_I$ equals}
$$
\gather
\sum_{\si}\ \frac{\sgn(\si)}{2^k\tss k\tss!}\cdot
\bigl(E_{i_{\si(1)},-i_{\si(2)}}-E_{i_{\si(2)},-i_{\si(1)}}\bigr)
\ldots
\bigl(E_{i_{\si(2k-1)},-i_{\si(2k)}}-E_{i_{\si(2k)},-i_{\si(2k-1)}}\bigr)
\\
=\ \sum_{\si}\ \frac{\sgn(\si)}{k\tss!}\cdot
E_{i_{\si(1)},-i_{\si(2)}}\ldots E_{i_{\si(2k-1)},-i_{\si(2k)}}
\endgather
$$
in $\UsoN\subset\UglN$.
Therefore the image $\ga(\Phi_I)\in\PD$ equals the sum
$$
\sum_{\si}\sum_{a_1\lc a_k=1}^m
\frac{\sgn(\si)}{k\tss!}\cdot
x_{a_1i_{\si(1)}}\d_{a_1,-i_{\si(2)}}\ldots\tss
x_{a_ki_{\si(2k-1)}}\d_{a_k,-i_{\si(2k)}}\,.
\Tag{4.5}
$$
Move each operator $\d_{a_p,-i_{\si(2p)}}$ in every
monomial of \(4.5) to the right, commuting it consecutively with the
multiplication operators 
$x_{a_{q}i_{\si(2q-1)}}$ for $q=p+1\lc k$.
Let us start with the operator $\d_{a_1,-i_{\si(2)}}$.
Take the commutator
$$
[\tss\d_{a_1,-i_{\si(2)}}\tss,\tss x_{a_2 i_{\si(3)}}\tss]=
\de_{a_1a_2}\,\de_{-i_{\si(2)},i_{\si(3)}}.
\nopagebreak
$$
The sum over $\si$ in (4.5) involves two commutators of this
form with opposite signs. Let $\bar\si\in S_{2k}$ is such that
$\bar\si(r)=\si(r)$ for $r\ne 2,3$ and
$\bar\si(2)=\si(3)$,$\bar\si(3)=\si(2)$.
Then $\sgn(\bar\si)=-\sgn(\si)$ whilst
$
[\tss\d_{a_1,-i_{\bar\si(2)}}\tss,\tss x_{a_2 i_{\bar\si(3)}}\tss]=
[\tss\d_{a_1,-i_{\si(2)}}\tss,\tss x_{a_2 i_{\si(3)}}\tss]\tss.
$
So by repeating this argument we get 
$$
\ga(\Phi_I)=
\sum_{\si}\sum_{a_1\lc a_k=1}^m
\frac{\sgn(\si)}{k\tss!}\cdot
x_{a_1i_{\si(1)}}\ldots x_{a_ki_{\si(2k-1)}}
\,
\d_{a_1,-i_{\si(2)}}\ldots\tss\d_{a_k,-i_{\si(2k)}}\,.
$$
This equality can be obviously rewritten as
$$
\ga(\Phi_I)\ =\ \sum_J
\sum_{a_1,\dots,a_k=1}^m
\frac{\sgn(JJ^\prime)}{k\tss!}\cdot
\det\tss[\tss x_{a_pj_q}\tss]\tss
\det\tss[\tss\d_{a_p,-j_q^{\tss\prime}}\tss]
$$
where $J$ and $J^\prime$ are the same sets as in \(4.4).
We now obtain the second relation in \(4.3) since  
both the determinants $\det\tss[\tss x_{a_pj_q}\tss]$
and $\det\tss[\tss\d_{a_p,-j_q^{\tss\prime}}\tss]$
are skew-symmetric with respect to permutations of the sequence
$(a_1\lc a_k)$.

Now let us show that the right hand side of the first relation
in \(4.3) belongs to the ring of $O_N$-invariants in $\UsoN$.
Consider the element of $\La^2(\CC^N)\ot\UsoN$
$$
\Phi\,=\,\sum_{i<j}\ (e_i\wedge e_j)\ot F_{i,-j}
\,=\,\frac12\tss\cdot\tss\sum_{ij}\ (e_i\wedge e_j)\ot F_{i,-j}\,.
$$
This element is invarant with respect to the diagonal action
of the group $O_N$. That follows from the invariance of the
element \(3.00). By the definition \(4.2) we have an equality
in the algebra $\La(\CC^N)\ot\UsoN$
$$
\frac{\,\Phi^k}{k\tss!}\,=\,\sum_I\ 
(e_{i_1}\wedge\ldots\wedge e_{i_{2k}})\ot\Phi_I\,.
$$
Extend the bilinear form $\langle\,\,,\,\rangle$ from
$\CC^N$ to the exterior power $\La^{2k}(\CC^N)$ by setting
$$
\langle\,
e_{i_1}\wedge\ldots\wedge e_{i_{2k}}
\,,\,
e_{-h_{2k}}\wedge\ldots\wedge e_{-h_1}\rangle
\,=\,
\de_{i_1h_1}\ldots\,\de_{i_{2k}h_{2k}}
$$
if 
$h_1<\ldots<h_{2k}$.
This extension is also $O_N$-invarant.
Extend it further to the $\CC$-bilinear form on $\La^{2k}(\CC^N)\ot\UsoN$
valued in $\UsoN$
$$
\langle\,\xi\ot X,\eta\ot Y\tss\rangle=
\langle\,\xi\,,\,\eta\tss\rangle\cdot XY\,,
$$
the latter form is $O_N$-equivariant.
This implies $O_N$-invariance of the sum in \(4.3)
$$
\sum_I\,\Phi_I\tss\Phi_\Is\,=\,
\langle\,\frac{\,\Phi^k}{k\tss!}\,,\,\frac{\,\Phi^k}{k\tss!}\,\rangle\,.
\Tag{4.6}
$$

Thus both sides of the first relation in \(4.3) belong to $\ZsoN$.
We will show that they have the same eigenvalue in any irreducible
$\soN$-module $V_\la$. By definition, eigenvalue of the left hand side
is the symmetric polynomial $(-1)^k\tss e_k(l_1^2\lc l_n^2\tss|\tss a)$
where we use the notation of Section 2 for $\g=\soN$.
Let us now determine the eigenvalue of the right hand side.
We will make use of the second equality in \(4.3).

Suppose that $m\ge n$, then every
$\soN$-module $V_\la$ is contained in $\P$ as a submodule. 
Indeed, there is a well known formula for a $\soN$-singular vector in
$\P$ of weight $\la=(\la_1\lc\la_n)$; see for instance
[\tss KV\,,\tss Section 6\tss]. For each $p=1\lc n$ denote
$$
\De_p=\det[\tss x_{ai}\tss]\,;
\qquad
a=1\lc\,p\,;
\quad
i=-n\lc\tss p\!-\!1\!-\!n\tss.
$$
One easily checks that the polynomial
$
v_\la=
\De_1^{\la_1-\la_2}\ldots\,\tss\De_{n-1}^{\la_{n-1}-\la_n}\tss\De_n^{\la_n}
\in\P
$
satisfies
$$
\align
F_{ij}\cdot v_\la&=0\,,
\,\ \ \quad\qquad\qquad\,i<j\,;
\\
F_{-p,-p}\cdot v_\la&=\la_{n-p+1}\,v_\la\,,
\quad 1\le p\le n
\endalign
$$
with respect to the action \(4.1).
Note that the polynomial $v_\la$ is of degree $|\tss\la\tss|$.
So due to the second relation in \(4.3) the right hand side of 
the first relation
vanishes in every module $V_\la$ with $|\la|<k$\tss;
see \(4.4).
Take the image of the left hand side of \(4.6)
with respect to the Harish-Chandra isomorphism.
Consider this image as a polynomial in $\la_1\lc\la_n$.
Due to Theorem 2.1 it now suffices to show that the leading component
of this polynomial is $e_k(\la_1^2\lc\la_n^2)$.

The leading terms of this polynomial arise from the summands of
$\Phi_I$ and $\Phi_\Is$ which depend only on the generators $F_{ii}$.
So here we may restrict the sum in \(4.6) to the sets $I$ of the form
$\{-p_k\lc\!-\!p_1,p_1\lc p_k\}$ where $1\le p_1<\ldots<p_k\le n$.
\text{Then we get $\Is=I.$\!} For any of these sets $I$ there
are $2^k\tss k\tss!$ summands in \(4.2) depending only on $F_{ii}$.
All these summands are equal to each other. \text{So the leading} terms
of our polynomial add up to
$$
\gather
\sum_{1\le p_1<\ldots<p_k\le n}\om\,
\bigl(F_{p_1p_1}^{\,2}\ldots F_{p_kp_k}^{\,2}\bigr)
=\sum_{1\le p_1<\ldots<p_k\le n}
\la_{n-p_1+1}^{\,2}\ldots\la_{n-p_k+1}^{\,2}
\\
=e_k(\la_1^2\lc\la_n^2)
\quad\square
\endgather
$$
\enddemo

\nt
The first equality in \(4.3) is an alternative
to Theorem 3.2 for $\g=\soN$.
However, for $\g=\soN$ there is apparently
no such alternative to Theorem 3.3.
When \text{$\g=\spN$} there is an analogue of
Theorem 4.1, but only for the element $D_k\in\ZspN$
as opposed to the element $C_k\in\ZspN$;
\text{see Section 5.}

\proclaim{Corollary 4.2}
Suppose that $N=2n$. Then we have an equality in $\UsoN$
$$
C_n=(-1)^n\cdot\bigl(\,\Pf\,[\tss F_{i,-j}\tss]\,\bigr)^2\,;
\qquad i,j=-n\lc\!-\!1,1\lc n.  
$$
\endproclaim

\nt
Let us now consider the operators from $\PD$ which commute with
the natural action of the group $O_N$ in $\P$. Description of this
commutant is a particular case of the theory of reductive dual pairs
due to R.\,Howe. Put $M=2m$. 
Let $\Esp_{\pm a,\pm b}$ be the matrix units in
$\EndCM$ corresponding to the standard basis of $e_{\pm a}$
in $\CC^{M}$. We will fix an embedding of the Lie algebra $\spM$
into $\glM$ similar to that used in Section 2 for $\spN$. So the
subalgebra $\spM\subset\glM$ will be spanned by the elements
$$
\gather
\Fsp_{ab}=-\Fsp_{-a,-b}=\Esp_{ab}-\Esp_{-b,-a}\,,
\\
\Fsp_{a,-b}=\Esp_{a,-b}+\Esp_{b,-a}\,,\ \quad
\Fsp_{-a,b}=\Esp_{-a,b}+\Esp_{-b,a}.
\endgather
\nopagebreak
$$
We also fix the Cartan decomposition of $\spM$ similar to that of $\spN$.
In particular, we will fix the Cartan subalgebra in $\spM$ with
the basis $\bigl(\tss\Fsp_{-m,-m}\,\lc\,\Fsp_{-1,-1}\bigr)$.
Any weight $\lap=(\la_1^\prime\lc\la_m^\prime)$ of that Cartan subalgebra
in $\spM$ will be taken with respect to this basis.
To distinguish the counterparts in $\ZspM$ of the elements
$C_1\lc C_n\in\ZsoN$ we will denote them here by
$C_1^{\tss\prime}\lc C_m^{\tss\prime}$.

One can define an action of the Lie algebra $\spM$ in the space $\P$ by
$$
\gather
\Fsp_{ab}\,\mapsto\,\sum_i\,x_{ai}\tss \d_{bi}+\frac{N}2\cdot\de_{ab}\,,
\\
\Fsp_{a,-b}\,\mapsto\,-\sum_i\,x_{ai}\tss x_{b,-i}\,,
\qquad
\Fsp_{-a,b}\,\mapsto\,\sum_i\,\d_{ai}\tss \d_{b,-i}\,.
\endgather
$$
Denote by $\gap$ the corresponding representation $\UspM\to\PD$.
The image of this representation coincides with the commutant in $\PD$
of the action of the group $O_N$ in $\P$ [\tss H\tss,\tss Theorem 8\tss].
In particular, we have the equality in $\PD$ 
$$ 
\ga\tss\bigl(\tss\ZsoN\bigr)=\gap\bigl(\tss\ZspM\bigr).
\Tag{4.7}
$$
We will describe the correspondence between the generators
\text{$C_1\lc C_n\in\ZsoN$}
and $C_1^{\tss\prime}\lc C_m^{\tss\prime}\in\ZspM$
implied by \(4.7); \text{cf. [\tss HU\tss,\tss Section 1\tss].}

As in Section 3 let $u$ be a complex variable.
Introduce the generating functions 
$$
\align
C(u)\,&=\,1\tss+\tss\sum_{k=1}^n\ 
\frac{C_k}
{\bigl(u^2-(\frac{N}2-1)^2\bigr)\ldots\bigl(u^2-(\frac{N}2-k)^2\bigr)}
\,\in\,\ZsoN\tss(u)\,,
\Tag{4.857}
\\
\hskip1cm
C^{\tss\prime}(u)\,&=\,1\tss+\tss\sum_{k=1}^m\ 
\frac{C_k^{\tss\prime}}
{\bigl(u^2-m^2\bigr)\ldots\bigl(u^2-(m-k+1)^2\bigr)}
\,\in\,\ZspM\tss(u)\,.
\endalign
$$
Also put in $\CC(u)$
$$
\al(u)\,=\,\prod_{a=1}^m\ \frac{u^2-(\frac{N}2-a)^2}{u^2-a^2}\,.
$$
Relations between the images
of \text{$C_1\lc C_n\in\ZsoN$}
and $C_1^{\tss\prime}\lc C_m^{\tss\prime}\in\ZspM$ in $\PD$ respectively 
under $\ga$ and $\gap$ are described by the following proposition.

\proclaim{Proposition 4.3}
We have an equality in $\PD\tss(u)$
$$
\al(u)\cdot\ga\bigl(C(u)\bigr)=\gap\bigl(C^{\tss\prime}(u)\bigr).
\Tag{4.75}
$$
\endproclaim

\demo{Proof}
The half-sum of positive roots for $\spM$ is $\rho^{\tss\prime}=(m\lc 1)$.
Let $l^{\tss\prime}_a=\la_a^\prime+\rho^{\tss\prime}_a$ be the analogues for
$\spM$ of the variables $l_p=\la_p+\rho_p$ for $\soN$. Recall that
here we have
$$
\rho=\tsize\bigl(\frac{N}2-1\lc\bigl\{\frac{N}2\bigr\}\bigr)
$$
where the braces mean the fractional part.
The images of $C(u)$ and
$C^{\tss\prime}(u)$~under the relevant Harish-Chandra isomorphisms are 
equal to
$$
\frac{
(u^2-l_1^{\tss2})
\ldots
(u^2-l_n^{\tss2})}
{(u^2-\rho_1^{\tss2})
\ldots
(u^2-\rho_n^{\tss2})
}
\ \quad\text{and}\ \quad
\frac{
\bigl(u^2-{l_1^{\tss\prime}}^2\bigr)
\ldots
\bigl(u^2-{l_m^{\tss\prime}}^{\!2}\bigr)
}{
\bigl(u^2-{\rho_1^{\tss\prime}}^2\bigr)
\ldots
\bigl(u^2-{\rho_m^{\tss\prime}}^{\!2}\bigr)
}
\Tag{4.8}
$$
respectively. By the definition of the elements $C_1\lc C_n$
the first equality follows from the identity 
\(1.6) applied to $t=u^2$ and
$(z_1\lc z_n)=(l_1^{\tss2}\lc l_n^{\tss2})$; see also 
\(2.88888888) and \(2.8). The second equality is obtained in a similar way.

Let us now employ the decomposition [\tss KV\tss,\tss Theorem 7.2\tss]
of the space $\PD$ into a direct sum of irreducible
$O_N\x\spM$-modules. By comparing the two rational functions in \(4.8)
one can now check that in each of these irreducible modules the eigenvalues of
$\al(u)\cdot C(u)\in\ZsoN\tss(u)$ and $C^{\tss\prime}(u)\in\ZspM\tss(u)$
coincide. This is a straightforward calculation, and we
refer to [\tss KV\tss,\tss Section 6\tss] for details 
\enddemos

\nt
By using the equality \(4.75) each of the images
$\gap(C_1^{\tss\prime})\lc\gap(C_m^{\tss\prime})$ in $\PD$
can be expressed as a linear combination of the operators
$\ga(C_1)\lc\ga(C_n)$; see the next theorem.
But for the latter operators we have already obtained
the explicit formulas \(4.3). Recall that for $k>n$ we have
$C_k=0$ by definition. Put $C_0=1$. 

\proclaim{Theorem 4.4}
Denote $d=m-\frac{N}2+1$. Then for any $k=1\lc m$ we have
$$
\gap(C^{\tss\prime}_k)=
\sum_{l=0}^k\ \ga\tss(C_l)\tss\tss\cdot
\prod_{s=l}^{k-1}\,
\frac{(m-s)\,(\frac{N}2-s)\,(d+l-s)}{(k-s)}\,\tss.
\Tag{4.9}
$$
\endproclaim 

\demo{Proof}
We will check that \(4.9) solves the equation \(4.75)
for $\gap(C_1^{\tss\prime})\lc\gap(C_m^{\tss\prime})$.
Let us subsitute \(4.9) into \(4.75) using the definitions of
the generating functions $C(u)$ and $C^{\tss\prime}(u)$.
Denote by $f_{kl}$ the product over $s=l\lc k-1$ at the
right hand side of \(4.9).
Then we have to check for any fixed $l=0,1\lc\min\tss(m,n)$
the equality 
$$
\tsize\prod_{r=1}^{m-l}\,\bigl(u^2-(r-d)^2\bigr)
\,=\,
\sum_{k=l}^m\,f_{kl}\,\cdot\prod_{r=1}^{m-k}\,(u^2-r^2)\,.
$$
Actually, the latter equality is valid for any $l=0,1\lc m$.
This can be easily verified by induction on $m-l=0,1,\,\ldots$
\enddemos

\nt
Observe that if $N=2n$ and $m>n$ then the coefficient $f_{kl}$
vanishes for any $k>n$. Therefore we have
the following corollary to Theorem 4.4.

\proclaim{Corollary 4.5}
Suppose $N=2n$ and $m>n$. Then
\text{$\gap(C_{n+1}^{\tss\prime})\lc\gap(C_m^{\tss\prime})=0$.} 
\endproclaim

\nt
Also observe that if $N=2n$ and $m\ge n-1$ then $f_{kl}=0$ for
$k>m-n+1+l$. In particular, we have another corollary;
\text{it also follows from Proposition 4.3 directly.}

\proclaim{Corollary 4.6}
Suppose that $N=2n$ and $m=n-1$. Then
$\gap(C_k^{\tss\prime})=\ga\tss(C_k)$ for any $k=1\lc m$.
\endproclaim

\nt
Note that due to Proposition 2.3 we have an equality in
$\ZsoN\tss[[\tss u\1\tss]]$
$$
1\tss+\tss\sum_{k=1}^\infty\ 
\frac{D_k}
{\bigl(u^2-(\frac{N}2)^2\bigr)\ldots\bigl(u^2-(\frac{N}2+k-1)^2\bigr)}
\,=C\,(u)\1
\,.
\Tag{5.*}
$$
Thus the series $C(u)\1$ can be regarded as a generating function
for the elements $D_1,D_2,\tss\ldots\in\ZsoN$; cf. the equality (5.11) 
for the Lie algebra $\g=\spN$ below.

Results of this section rely on the fact that
real Lie groups \text{$O_N(\RR)$ and $Sp_{\tss M}(\RR)$} form a
reductive dual pair inside $Sp_{\tss MN}(\RR)$; see [\tss EHW\tss].
The case of another dual 
\text{pair of $Sp_{\tss N}(\RR)$ and $O^{\,\ast}_M(\RR)$
in $Sp_{\tss MN}(\RR)$ will be considered in the next section.}
 

\kern15pt
\line{\bf5.\ Capelli identity for the Lie algebra $\spN$\hfill}
\section{\,}\kern-20pt

\nt
In this section we will keep fixed two positive integers $m$ and $N$.
We will consider a natural representation of the enveloping algebra
$\UspN$ by differential operators with polynomial coefficients on the
space $\CMN$. Results of this section will be similar to those of
Section 4. But here for any $k=1,2,\ldots$ we will give
an explicit formula for the differential
operator corresponding to the element $D_k\in\ZspN$ rather than 
$C_k\in\ZspN$. This formula may be still regarded as an analogue for $\spN$
of the Capelli identity [\tss C\tss] for the Lie algebra $\glN$;
cf. [\tss N1\tss,\tss Example 2\tss].

Put $N=2n$. Let the indices $i,j$ run
through the set $\{\tss-\tss n\lc\!-\!1,1\lc n\tss\}$.
The indices $a,b$ will again run through the set $\{\tss1\lc m\tss\}$.
Let $\P$ be the same ring of polynomials in the variableas $x_{ai}$ 
as introduced in the previous section.
We will use the embedding of the group $Sp_N$ into $GL_N$ chosen in
Section 2. The restriction to $\spN$ of the action \(4.0) of $\glN$ in $\P$
then takes the form
$$
F_{ij}\,\mapsto\,\sum_a\ \bigl(\,x_{ai}\d_{aj}-
\sgn(ij)\cdot
x_{a,-j}\d_{a,-i}\tss\bigr)\,.
\Tag{5.1}
$$  
Thus we get a representation $\UspN\to\PD$ which we will again denote by
$\ga$. For each $k=1,2,\ldots$ we will give an explicit formula
for the operator
$\ga\tss(D_k)\in\PD$.

Put $\Ft_{ij}=\sgn(i)\cdot F_{ij}$. Then 
by definition we have $\Ft_{i,-j}=\Ft_{j,-i}$.
Let $I$ be any sequence of indices $i$ with $2k$ \text{elements.}
Suppose these elements are $i_1\le\ldots\le i_{2k}$.
\text{Take an element of the algebra $\UspN$} 
$$
\Psi_I=\sum_{\si}\ \frac1{2^k\tss k\tss!}\cdot
\Ft_{i_{\si(1)},-i_{\si(2)}}\ldots\Ft_{i_{\si(2k-1)},-i_{\si(2k)}}
\Tag{5.2}
$$
where $\si$ runs through the symmetric group $S_{2k}$.
The element \(5.2) may be called the {\it Hafnian} of the
symmetric matrix 
$[\tss\Ft_{i_p,-i_q}\tss]$ with the indices $p,q=1\lc 2k$.
This term was devised by E.\,Caianiello for a symmetric matrix
with commuting entries; see [\tss K\tss,\tss Section 2\tss].
The factors of the monomial in the sum \(5.2) do not commute
in general. Still we will use this term and write
$$
\Psi_I=\Hf\,[\tss\Ft_{i_p,-i_q}\tss]\,;
\ \quad p,q=1\lc 2k.  
$$

Let us denote $\Is=(-\tss i_{2k}\lc\!-\!i_1)$. 
Further, let $A$ be any collection of indices $a$ with $k$ elements.
Suppose these elements are $a_1\lc a_k$. Denote in $\PD$  
$$
\Theta_{AI}\,=\,
\sum_J\ 
\sgn\tss(j_1\ldots j_k)\cdot
\per\tss[\tss x_{a_pj_q}\tss]\tss
\per\tss[\tss\d_{a_p,-j_q^{\tss\prime}}\tss]
\Tag{5.4}
$$
where the sum is taken over all partitions of $I$ into subsequences
$J=(\tss j_1\lc j_k\tss)$ and
$J^\prime=(\tss j_1^{\tss\prime}\lc j_k^{\tss\prime}\tss)$
each of length $k$. So $j_1\le\ldots\le j_k$
\text{and $j_1^{\tss\prime}\le\ldots\le j_k^{\tss\prime}$.} 
Here the permanents are taken with respect to the indices $p,q=1\lc k$.
Let $f_{\pm1}\lc f_{\pm n}$ be the 
multiplicities of the numbers \text{$\pm1\lc\pm\!n$ in $I$}
respectively. Let $d_1\lc d_m$ be the
multiplicities of $1\lc m$ in $A$.
The next statement is an analogue of Theorem 4.1, it
provides an explicit formula for $\ga(D_k)\in\PD$.
\nopagebreak
\proclaim{Theorem 5.1}
We have the equalities respectively in $\UspN$ and $\PD$
$$
D_k=(-1)^k\cdot\sum_I\,
\frac{\,\sgn(i_1\ldots i_{2k})\cdot\Psi_I\tss\Psi_\Is}
{f_1\tss!\tss f_{-1}\tss!\,\tss\ldots\,f_n\tss!\tss f_{-n}\tss!}
\ \ \text{and}\ \ \,
\ga(\Psi_I)=\sum_A\,
\frac{\Theta_{AI}}{d_1\tss!\tss\ldots\,d_m\tss!}\,.
\Tag{5.3}
$$
\endproclaim
 
\demo{Proof}
We start with verifying the second equality in \(5.3).
\text{By definition $\!\Psi_I\!$ equals}
$$
\align
\frac1{2^k\tss k\tss!}\,\cdot\sum_{\si}\ 
&\bigl(\tss
\sgn(i_{\si(1)})\cdot E_{i_{\si(1)},-i_{\si(2)}}+
\sgn(i_{\si(2)})\cdot E_{i_{\si(2)},-i_{\si(1)}}
\bigr)
\x\ldots\x
\\
&\bigl(\tss
\sgn(i_{\si(2k-1)})\cdot E_{i_{\si(2k-1)},-i_{\si(2k)}}+
\sgn(i_{\si(2k)})\cdot E_{i_{\si(2k)},-i_{\si(2k-1)}}
\bigr)
\,=
\\
\frac1{k\tss!}\,\cdot\sum_{\si}\,
&\ \sgn(i_{\si(1)}\tss i_{\si(3)}\ldots i_{\si(2k-1)})\cdot
E_{i_{\si(1)},-i_{\si(2)}}\ldots E_{i_{\si(2k-1)},-i_{\si(2k)}}
\endalign
$$
in $\UspN\subset\UglN$.
Therefore the image $\ga(\Psi_I)\in\PD$ equals the sum
$$
\sum_{\si}\sum_{a_1\lc a_k=1}^m
\frac{\,\sgn(i_{\si(1)}\ldots i_{\si(2k-1)})}{k\tss!}\cdot
x_{a_1i_{\si(1)}}\d_{a_1,-i_{\si(2)}}\ldots\tss
x_{a_ki_{\si(2k-1)}}\d_{a_k,-i_{\si(2k)}}\,.
$$
Move to the right each operator $\d_{a_p,-i_{\si(2p)}}$ in every
monomial above, commuting it consecutively with the
multiplication operators 
$x_{a_{q}i_{\si(2q-1)}}$ for $q=p+1\lc k$.
Let us start with the operator $\d_{a_1,-i_{\si(2)}}$.
Take the commutator
$$
[\tss\d_{a_1,-i_{\si(2)}}\tss,\tss x_{a_2 i_{\si(3)}}\tss]=
\de_{a_1a_2}\,\de_{-i_{\si(2)},i_{\si(3)}}.
$$
The above sum over $\si$ contains two such commutators
with opposite signs. Indeed, if $\bar\si\in S_{2k}$ is such that
$\bar\si(r)=\si(r)$ for $r\ne 2,3$ and
$
\bar\si(2)=\si(3),\bar\si(3)=\si(2)
$
then 
$$
[\tss\d_{a_1,-i_{\bar\si(2)}}\tss,\tss x_{a_2 i_{\bar\si(3)}}\tss]=
[\tss\d_{a_1,-i_{\si(2)}}\tss,\tss x_{a_2 i_{\si(3)}}\tss]
$$
whilst 
$$
\sgn(i_{\bar\si(1)}\tss i_{\bar\si(3)}\tss\ldots\tss i_{\bar\si(2k-1)})=
-\sgn(i_{\si(1)}\tss i_{\si(3)}\tss\ldots\tss i_{\si(2k-1)})
$$
for $i_{\si(2)}=-i_{\si(3)}$. 
So by repeating this argument we bring $\ga(\Psi_I)$ to the form
$$
\sum_{\si}\sum_{a_1\lc a_k=1}^m
\frac{\sgn(i_{\si(1)}\ldots i_{\si(2k-1)})}{k\tss!}\cdot
x_{a_1i_{\si(1)}}\ldots x_{a_ki_{\si(2k-1)}}
\,
\d_{a_1,-i_{\si(2)}}\ldots\tss\d_{a_k,-i_{\si(2k)}}\,.
$$
The latter sum in $\PD$ can be obviously rewritten as
$$
\ga(\Psi_I)\ =\ \sum_J
\sum_{a_1,\dots,a_k=1}^m
\frac{\sgn(j_1\ldots j_k)}{k\tss!}\cdot
\per\tss[\tss x_{a_pj_q}\tss]\tss
\per\tss[\tss\d_{a_p,-j_q^{\tss\prime}}\tss]
$$
where $J$ and $J^\prime$ are the same as in \(5.4).
We now obtain the second relation in \(5.3) since  
both the permanents $\per\tss[\tss x_{a_pj_q}\tss]$
and $\per\tss[\tss\d_{a_p,-j_q^{\tss\prime}}\tss]$
are symmetric with respect to permutations of the sequence
$(a_1\lc a_k)$.

Now let us show that the right hand side of the first relation
in \(5.3) belongs to the centre of $\UspN$.
Consider the element of $\operatorname{S}^2(\CC^N)\ot\UspN$
$$
\Psi\,=\,\frac12\tss\cdot\tss\sum_{ij}\ e_i\tss e_j\ot\Ft_{i,-j}\,.
\nopagebreak
$$
This element is invarant with respect to the diagonal action
of the group $Sp_N$. That again follows from the invariance of the
element \(3.00). By the definition \(5.2) we have an equality
in the algebra $\operatorname{S}\tss(\CC^N)\ot\UsoN$
$$
\frac{\,\Psi^k}{k\tss!}\,=\,\sum_I\ 
\frac{e_{i_1}\ldots e_{i_{2k}}}
{f_1\tss!\tss f_{-1}\tss!\,\tss\ldots\,f_n\tss!\tss f_{-n}\tss!}
\ot\Psi_I\,.
$$
\text{Extend the bilinear form $\langle\,\,,\,\rangle$ from $\CC^N$ to the
symmetric power $\operatorname{S}^{2k}\tss(\CC^N)$ by}
$$
\langle\,
e_{i_1}\ldots e_{i_{2k}}
\tss,\tss
e_{-h_{2k}}\ldots e_{-h_1}\rangle
\tss=\tss
\de_{i_1h_1}\ldots\,\de_{i_{2k}h_{2k}}\cdot
{f_1!f_{-1}!\tss\ldots f_n!f_{-n}!}
\cdot\sgn(i_1\ldots i_{2k})\,.
$$
This extension is also $Sp_N$-invarant.
Extend it further to the $\CC$-bilinear form on
$\operatorname{S}^{2k}(\CC^N)\ot\UspN$
valued in $\UspN$
$$
\langle\,\xi\ot X\tss,\tss\eta\ot Y\tss\rangle=
\langle\,\xi\tss,\tss\eta\tss\rangle\cdot XY\,,
$$
the latter form is $Sp_N$-equivariant.
\text{It implies $Sp_N$-invariance of the sum in \(5.3)}
$$
\sum_I\ 
\frac{\,\sgn(i_1\ldots i_{2k})\cdot\Psi_I\tss\Psi_\Is}
{f_1\tss!\tss f_{-1}\tss!\,\tss\ldots\,f_n\tss!\tss f_{-n}\tss!}
\,=\,
\langle\,\frac{\,\Psi^k}{k\tss!}\,,\tss\frac{\,\Psi^k}{k\tss!}\,\rangle\,.
\Tag{5.6}
$$

Thus both sides of the first relation in \(5.3) belong to $\ZspN$.
We will show that they have the same eigenvalue in any irreducible
$\spN$-module $V_\la$. By definition, the eigenvalue of the left hand side
is the symmetric polynomial $h_k(l_1^2\lc l_n^2\tss|\tss a)$
where we use the notation of Section 2 for $\g=\spN$.
Let us now determine the eigenvalue of the right hand side.
We will make use of the second equality in \(5.3).

Suppose that $m\ge n$, then every
$\spN$-module $V_\la$ is contained in $\P$ as a submodule. 
Indeed, the vector $v_\la\in\P$ defined in the proof
of Theorem 4.1 is also singular of weight $\la=(\la_1\lc\la_n)$
\text{with respect to the action \(5.1) of $\spN$. That is,}
$$
\align
F_{ij}\cdot v_\la&=0\,,
\,\ \ \quad\qquad\qquad\,i<j\,;
\\
F_{-p,-p}\cdot v_\la&=\la_{n-p+1}\,v_\la\,,
\quad 1\le p\le n\,
\endalign
$$
again. The polynomial $v_\la$ is of degree $|\tss\la\tss|$.
By second relation in \(5.3) the right hand side of the first relation
vanishes in every module $V_\la$ with $|\tss\la\tss|<k$\tss;
\text{see \(5.4).}

Take the image of the left hand side of \(5.6)
with respect to the Harish-Chandra isomorphism. 
Let us consider this image as a polynomial in $\la_1\lc\la_n$.
Due to Theorem 2.1 it now suffices to show that the leading component
of this polynomial is $(-1)^k\cdot h_k(\la_1^2\lc\la_n^2)$.

The leading terms of this polynomial arise from the summands of
$\Psi_I$ and $\Psi_\Is$ which depend only on the generators $F_{ii}\in\spN$.
So here we may restrict the sum  in \(5.6) to the sequences $I$ of the form
$\{-p_k\lc\!-\!p_1,p_1\lc p_k\}$ where $1\le p_1\le\ldots\le p_k\le n$.
Then we get $\Is=I$.
For any of these sequences $I$ there
are $2^{\tss k}\tss k\tss!\,\cdot f_1\tss!\ldots f_n\tss!$ summands in \(5.2)
depending only on $F_{ii}$.
All these summands are equal to each other.
So the leading terms of our polynomial \text{add up to}
$$
\gather
\sum_{1\le p_1\le\ldots\le p_k\le n}\om\,
\bigl(F_{p_1p_1}^{\,2}\ldots F_{p_kp_k}^{\,2}\bigr)\cdot(-1)^k
=\!\!\sum_{1\le p_1\le\ldots\le p_k\le n}
\la_{n-p_1+1}^{\,2}\ldots\la_{n-p_k+1}^{\,2}\cdot(-1)^k
\\
=(-1)^k\tss h_k(\la_1^2\lc\la_n^2)
\quad\square
\endgather
$$
\enddemo

\nt
The first equality in \(5.3) is an alternative
to Theorem 3.3 for $\g=\spN$.
However, when $\g=\spN$ there is apparently 
no such alternative to Theorem 3.2.

Let us now consider the operators from $\PD$ which commute with
the action \(5.1) of the Lie algebra $\spN$ in $\P$. Description of this
commutant is another particular case of the theory of reductive dual pairs
[\tss H\tss].
Put $M=2m$. As in Section~4,  
denote by $\Esp_{\pm a,\pm b}$ the matrix units in
$\EndCM$ corresponding to the standard basis of $e_{\pm a}$
in $\CC^{M}$. Fix an embedding of the Lie algebra $\soM$
into $\glM$ similar to that used in Section 2 for $\soN$.
\text{So the subalgebra $\soM\subset\glM$ will be spanned by}
$$
\gather
\Fsp_{ab}=-\Fsp_{-a,-b}=\Esp_{ab}-\Esp_{-b,-a}\,,
\\
\Fsp_{a,-b}=\Esp_{a,-b}-\Esp_{b,-a}\,,
\ \quad\Fsp_{-a,b}=\Esp_{-a,b}-\Esp_{-b,a}.
\endgather
$$
We also fix the Cartan decomposition of $\soM$ similar to that of $\soN$.
In particular, we will fix the Cartan subalgebra in $\soM$ with
the basis $\bigl(\tss\Fsp_{-m,-m}\,\lc\,\Fsp_{-1,-1}\bigr)$.
Any weight $\lap=(\la_1^\prime\lc\la_m^\prime)$ of that Cartan subalgebra
in $\soM$ will be taken with respect to this basis.
To distinguish the counterparts in $\ZsoM$ of the elements
$D_1,D_2\tss,\tss\ldots\in\ZspN$ we will denote them by
$D_1^{\tss\prime},D_2^{\tss\prime}\tss,\tss\ldots$ in this section.

One can define an action of the Lie algebra $\soM$ in the space $\P$ by
$$
\gather
\Fsp_{ab}\,\mapsto\,\sum_i\,x_{ai}\tss \d_{bi}+n\cdot\de_{ab}\,,
\\
\Fsp_{a,-b}\,\mapsto\,\sum_i\,\sgn(i)\cdot x_{ai}\tss x_{b,-i}\,,
\qquad
\Fsp_{-a,b}\,\mapsto\,\sum_i\,\sgn(i)\cdot\d_{ai}\tss \d_{b,-i}\,.
\endgather
$$
Denote by $\gap$ the corresponding representation $\UsoM\to\PD$.
The image of this representation coincides with the commutant in $\PD$
of the action of the Lie algebra $\spN$ in $\P$
[\tss H\tss,\tss Theorem 8\tss].
In particular, we have the equality in $\PD$ 
$$ 
\ga\tss\bigl(\tss\ZspN\bigr)=\gap\bigl(\tss\ZsoM\bigr).
\Tag{5.7}
$$
We will describe the correspondence between the generators
\text{$D_1,D_2\tss,\tss\ldots\in\ZspN$}
and $D_1^{\tss\prime},D_2^{\tss\prime}\tss,\tss\ldots\in\ZsoM$
implied by \(5.7).
 
As in Section 4 let $u$ be a complex variable.
Introduce the generating functions 
$$
\align
D(u)\,&=\,1\tss+\tss\sum_{k=1}^\infty\ 
\frac{D_k}
{\bigl(u^2-(n+1)^2\bigr)\ldots\bigl(u^2-(n+k)^2\bigr)}
\,\in\,\ZspN\tss[[\tss u\1\tss]]\,,
\\
D^{\tss\prime}(u)\,&=\,1\tss+\tss\sum_{k=1}^\infty\ 
\frac{D_k^{\tss\prime}}
{\bigl(u^2-m^2\bigr)\ldots\bigl(u^2-(m+k-1)^2\bigr)}
\,\in\,\ZsoM\tss[[\tss u\1\tss]]\,.
\endalign
$$
Also put in $\CC(u)$
$$
\be(u)\,=\,\prod_{a=1}^m\ \frac{u^2-(a-1)^2}{u^2-(n-a+1)^2}\,.
$$
Relations between the images
of \text{$D_1,D_2,\tss\ldots\in\ZspN$}
and $D_1^{\tss\prime},D_2^{\tss\prime},\tss\ldots\in\ZsoM$
in $\PD$ respectively 
under $\ga$ and $\gap$ can be described as follows.
\nopagebreak
\proclaim{Proposition 5.2}
We have an equality in $\PD\tss(u)$
$$
\be(u)\cdot\ga\bigl(D(u)\bigr)=\gap\bigl(D^{\tss\prime}(u)\bigr).
\Tag{5.75}
$$
\endproclaim

\demo{Proof}
It will be similar to that of Proposition 4.3.
The half-sum of positive roots for $\soM$ is $\rho^{\tss\prime}=(m-1\lc 0)$.
Let $l^{\tss\prime}_a=\la_a^\prime+\rho^{\tss\prime}_a$ be the analogues for
$\soM$ of the variables $l_p=\la_p+\rho_p$ for $\spN$. Recall that
here we have $\rho=(n\lc 1)$.
The images of $D(u)$ and
$D^{\tss\prime}(u)$~under the relevant Harish-Chandra \text{isomorphisms are} 
$$
\frac{
(u^2-\rho_1^{\tss2})
\ldots
(u^2-\rho_n^{\tss2})
}{
(u^2-l_1^{\tss2})
\ldots
(u^2-l_n^{\tss2})
}
\ \quad\text{and}\ \quad
\frac{
\bigl(u^2-{\rho_1^{\tss\prime}}^2\bigr)
\ldots
\bigl(u^2-{\rho_m^{\tss\prime}}^{\!2}\bigr)
}{
\bigl(u^2-{l_1^{\tss\prime}}^2\bigr)
\ldots
\bigl(u^2-{l_m^{\tss\prime}}^{\!2}\bigr)
}
\Tag{5.8}
$$
respectively. By the definitions of $D_1,D_2,\tss\ldots$
here the first equality follows from
\(1.7) applied to $t=u^2$ and
$(z_1\lc z_n)=(l_1^{\tss2}\lc l_n^{\tss2})$; see also 
\(2.88888888) and \(2.8).
The second equality is obtained in a similar way.

Let us now employ
the decomposition [\tss EHW\tss,\tss Section 9\tss]
of the space $\PD$ into a direct sum of irreducible
$\spN\x\soM$-modules. By comparing the two rational functions in \(5.8) one
can now check that in each of these irreducible modules the eigenvalues of
$\be(u)\cdot D(u)\in\ZspN\tss(u)$ and $D^{\tss\prime}(u)\in\ZsoM\tss(u)$
coincide. This is again a straightforward calculation, and we will omit
the details 
\enddemos

\nt
By using the equality \(5.75) each of the images
$\gap(D_1^{\tss\prime})\tss,\gap(D_2^{\tss\prime})\tss,\tss\ldots$ in $\PD$
can be expressed as a linear combination of the operators
$\ga(D_1)\tss,\ga(D_2)\tss,\tss\ldots\,$; see the next theorem.
For the latter operators we already have
the formulas \(5.3). Put $D_0=1$.

\proclaim{Theorem 5.3}
Put $d=n-m+1$. Then for any $k=1,2,\ldots$ we have
$$
\gap(D^{\tss\prime}_k)=
\sum_{l=0}^k\ \ga\tss(D_l)\tss\cdot
\prod_{s=l}^{k-1}\,
\frac{(m+s)\,(n+s)\,(d-k+s+1)}{(s-l+1)}\,\tss.
\Tag{5.9}
$$
\endproclaim 

\demo{Proof}
By Proposition 5.2 and the definitions of $D(u),D^{\tss\prime}(u)$
the differential operator $\gap(D_k^{\tss\prime})$ is a linear
combination of $\ga(D_1)\lc\ga(D_k)$ where the coefficients
are certain rational functions in $m,n$. Therefore it suffices
to prove \(5.9) only for $d\ge0$. Suppose this is the case.
We will check that \(5.9) solves the equation \(5.75) for 
$\gap(D_1^{\tss\prime})\tss,\gap(D_2^{\tss\prime})\tss,\tss\ldots\,\tss$.
Let us substitute \(5.9) into \(5.75) using the definitions of
the generating functions $D(u)$ and $D^{\tss\prime}(u)$.
Denote by $g_{kl}$ the product over $s=l\lc k-1$ at the
right hand side of \(5.9). Observe that $g_{kl}=0$ for
any $k>d+l$. So we have to check for any fixed $l=0,1,2,\ldots$
the equality in $\CC(u)$
$$
\prod_{r=0}^{m+l-1}\,\frac1{u^2-(r+d)^2}\,=\,
\sum_{k=l}^{l+d}\ g_{kl}\,\cdot\!\!
\prod_{r=0}^{m+k-1}\!\frac1{u^2-r^2}
$$
Actually, the latter equality is valid for any $l=-m,1-m,\ldots\,$.
This can be easily verified by induction on $m+l=0,1,\,\ldots$
\enddemos

\nt
We also have a direct corollary to Proposition 5.2, it is an analogue of
\text{Corollary 4.6.} 

\proclaim{Corollary 5.4}
Suppose $n=m-1$. Then
$\gap(D_k^{\tss\prime})=\ga\tss(D_k)$ \text{for any $k=1,2,\ldots$.}
\endproclaim

\nt
Note that due to Proposition 2.3 we have an equality in
$\ZspN\tss[[\tss u\1\tss]]$
$$
1\tss+\tss\sum_{k=1}^n\ 
\frac{C_k}
{\bigl(u^2-n^2\bigr)\ldots\bigl(u^2-(n-k+1)^2\bigr)}
\,=D\,(u)\1
\,.
\Tag{5.857}
$$
Similarly to \(4.857), denote this generating function by $C(u)$.
Then by definition we have the equality in $\ZspN\tss[[\tss u\1\tss]]$
$$
1\tss+\tss\sum_{k=1}^\infty\ 
\frac{D_k}
{\bigl(u^2-(n+1)^2\bigr)\ldots\bigl(u^2-(n+k)^2\bigr)}
\,=C\,(u)\1
\,.
\Tag{4.*}
$$
In the next section we will give another formula for the generating
function $C(u)$.


\kern15pt
\line{\bf6.\ Quantum determinant of the element $F(u)$\hfill}
\section{\,}
\kern-20pt

\nt
Here we will give an explicit formula for the generating
function $C(u)\in\Zg(u)$.
This function was defined by \(4.857) for
$\g=\soN$ and by \(5.857) for $\g=\sp_N$. We will use the notation and some
results of Section 3. In that section we introduced in particular
the rational function of $u$ valued in the algebra $\EndCN^{\ot N}\ot\Ug$
$$
F_{(1^N)}(u)=(A_N\ot1)\cdot\,\prod_{q=1}^N\ 
\biggl(
\Bigl(\tss
1+\frac{Q_{1q}+\ldots+Q_{q-1,q}}{2u-q+1}
\tss\Bigr)\ot1\cdot F_q(u-q+1)
\biggr)
\nopagebreak
$$
where the factors corresponding to the indices $q=1\lc N$ are arranged
from left to right. Recall that the element $F(u)\in\EndCN\ot\Ug\tss(u)$
\text{is given by \(3.02).}

Put $\ep(u)=(2u+1)/(2u-N+1)$ for $\g=\spN$ and $\ep(u)=1$ for $\g=\soN$.
By definition $F_{(1^N)}(u)$ is divisible by $A_N\ot1$ on the left.
Due to Proposition 3.7 and to the 
decomposition \(3.04) the element $F_{(1^N)}(u)$
is also divisible by $A_N\ot1$ on the right. To show this one should use
the relation \(3.111) repeatedly. But the projector $A_N$ is one-dimensional.
Therefore
$$
F_{(1^N)}(u)\,=\,\ep(u)\cdot A_N\ot\Cb(u)
\Tag{6.2}
$$
for certain rational function $\Cb(u)$ valued in $\Ug$.
Moreover,  due to Proposition~3.4 we have $\Cb(u)\in\Zg\tss(u)$.
Note that the definition of $F_{(1^N)}(u)$
yields an explicit formula for $\Cb(u)$.
The function $\Cb(u)$ is called the {\it quantum determinant} of
the element $F(u)$, or the {\it Sklyanin determinant\tss};
see [O] and [\tss MNO\tss,\tss Section 4\tss].
The reason for incorporating the factor $\ep(u)\in\CC(u)$
to the definition \(6.2)
of $\Cb(u)$ will become clear in the proof of Theorem 6.2.

That proof will be based on the following proposition.
Consider the polynomial in $u$ valued in the algebra $\EndCN^{\ot N}\ot\UglN$
$$
(A_N\ot1)\cdot E_1(u)\,E_2(u-1)\,\ldots\,E_N(u-N+1)\tss.
\nopagebreak
\Tag{6.61}
$$ 
Using the decomposition \(3.04) and the relation \(3.81) the product
\(6.61) can be rewritten as
$$
E_N(u-N+1)\,\ldots\,E_2(u-1)\,E_1(u)\cdot(A_N\ot1)\tss.
$$
Therefore \(6.61) equals $A_N\ot H(u)$ for a certain polynomial $H(u)$
valued in $\UglN$. This polynomial is called the {\it quantum determinant} of
the element $E(u)$. For comments see [\tss MNO\tss,\tss Section 2\tss] 
and references
therein. Arguments similar to those used in the proof of Proposition 3.4
show that actually $H(u)\in\ZglN\tss(u)$. 

\nt
Further, similarly to the proof of Corollary 3.11 one demonstrates that
$$
(A_N\ot1)\cdot\Et_1(u-N+1)\,\ldots\,\Et_{N-1}(u-1)\,\Et_N(u)=
A_N\ot H(u)\tss.
\Tag{6.62}
$$

Let $\nu$ be any partition
with $\ell(\nu)\le N$. Consider the corresponding irreducible
$\glN$-module $U_\nu$. The following proposition is well known,
see for instance [\tss H\tss].

\proclaim{Proposition 6.1}
The eigenvalue of $H(u)$ in the module $U_\nu$ is
$$
\prod_{q=1}^N\ (\tss\nu_q+N-q-u\tss)\tss.
$$
\endproclaim

\demo{Proof}
Let us fix the Borel subalgebra in the Lie algebra $\glN$ 
generated by the elements $E_{ij}$ with $i\le j$. Further,
fix the basis 
$
\bigl(E_{-n,-n}\,\lc\,E_{nn}\bigr)
$
in the Cartan subalgebra of $\glN$.
Then the module $U_\nu$ is of the highest weight $(\nu_1\lc \nu_N)$.
Let $\xi\in U_\nu$ be a highest weight vector. By considering the action
of the element \(6.61) on the vector 
$$
e_n\ot\ldots\ot e_{-n}\ot\xi\tss\in\tss\EndCN^{\ot N}\ot U_\nu
$$
and by using the definition of $H(u)\in\ZglN\tss(u)$
we get the required statement
\enddemos
 
\nt
The following theorem is the main result of the present section; cf. [M1]
and [NO]. 

\proclaim{Theorem 6.2}
We have an equality in the algebra $\Zg\tss(u)$
$$
\tsize
C(u)\,=\,\Cb(u+\frac{N}2-\frac12)\tss\cdot
\dsize
\prod_{q=1}^N\,
\tsize
\bigl(\tss\frac{N}2+\frac12-q-u-\eta\tss\bigr)\1
\,.
\Tag{6.3}
$$
\endproclaim

\demo{Proof}
Let $V_\la$ be an irreducible $\g$-module of highest weight
$\la=(\la_1\lc \la_n)$. As in Section 2 we assume that
$\la_n\ge0$ so that $\la$ is a partition with $\ell(\la)\le n$.
Due to \(1.6) and \(2.8) the eigenvalue of the element $C(u)\in\Zg(u)$
in $V_\la$ equals 
$$
\prod_{p=1}^n\ 
\frac{\tss u^2-(\la_p+\rho_p)^2}{u^2-\rho_p^{\tss2}}\,.
\Tag{6.4}
$$
We have to show that the element of $\Zg(u)$ at the right hand
side of \(6.3) has in $V_\la$ the same eigenvalue \(6.4). As well as in
Section 3, fix any embedding of the $\glN$-module $U_\la$ to 
$(\CC^N)^{\ot\tss l}$ where $l=|\tss\la\tss|$.
Let $V\subset(\CC^N)^{\ot\tss l}$ be the subspace of traceless tensors.
Recall that $U_\la\cap V=V_\la$ unless $\g=\so_{2n}$ and $\la_n\neq0$, 
in the latter case $U_\la\cap V=V_\la\oplus V_\las$. As well as in Section 3,
by using Lemma 3.8 we can show that the restriction of the action of
$
F_{(1^N)}(u)\in\EndCN^{\ot N}\ot\Ug(u)
$
to $(\CC^N)^{\ot\tss l}\ot V$ coincides with that of
the element from $\EndCN^{\ot N}\ot\UglN(u)$
$$
\gather
(A_N\ot1)\,\cdot\,\prod_{q=1}^{2k}\,\tss\Et_q(\eta-u+q-1)\,\,\x 
\Tag{6.1313}
\\
\prod_{q=1}^N\,
\biggr(\!
\Bigl(\tss1+\frac{Q_{1q}+\ldots+Q_{q-1,q}}{2u-q+1}\tss\Bigr)\ot1
\tss\!\biggr)
\,\cdot\,
\prod_{q=1}^N\,\frac1{u-q+1-\eta}\,\,\x
\\
(A_N\ot1)\,\cdot\,\prod_{q=1}^N\,\tss E_q(\eta+u-q+1)\quad  
\endgather
\nopagebreak
$$
where the non-commuting
factors corresponding to the indices $q=1\lc N$ are as

\nt
usual arranged from left to right. But the element \(6.1313) equals
$$
\biggl(
A_N\cdot
\prod_{q=1}^N\,
\Bigl(\tss1+\frac{Q_{1q}+\ldots+Q_{q-1,q}}{2u-q+1}\tss\Bigr)
\!\biggr)
\tss\ot\tss
\frac
{\tss H(\eta-u+N-1)\tss H(\eta+u)\tss}
{(u-\eta)\tss\ldots\tss(u-N+1-\eta)}
\,;
\Tag{6.63}
$$
see \(6.62). Further, the first tensor factor in \(6.63)
equals $\ep(u)\cdot A_N$. This can be verified by direct calculation; it
has been performed in [\tss MNO\tss,\tss Section 4\tss]. So by the
definition \(6.2) the eigenvalue of $\Cb(u)$ in the $\g$-module
$V_\la$ coincides with the eigenvalue in the $\glN$-module $U_\la$
of the second tensor factor in \(6.63). Due to Proposition 6.1 the latter
eigenvalue equals
$$
\prod_{p=1}^n\
\frac{(\la_p-\eta-p+u+1)\,(\la_p-\eta+N-p-u)}{u-\eta-p+1}
\,\ \cdot
\prod_{q=n+1}^N(N-q-\eta-u)\,.
$$
Since $\rho_p=\frac{N}2+\frac12-\eta-p$ for
$\g=\so_{2n}\,,\,\so_{2n+1}\,,\,\sp_{2n}$ we obtain
that $\Cb(u+\frac{N}2-\frac12)$ takes in $V_\la$ the eigenvalue
$$
\prod_{p=1}^n\,
\frac{\tss u^2-(\la_p+\rho_p)^2}{u^2-\rho_p^{\tss2}}
\,\ \cdot\ 
\prod_{q=1}^N\,(N-q-\eta-u)
\quad\square
$$
\enddemo

\nt
This theorem was obtained by the second author in 1991.
Another proof of this theorem was given by the first author in [M1].
It uses a certain map $S_N\to S_{N-1}$ and 
provides one more formula for the generating function
\text{$C(u)\in\Zg(u)$.}
\enddemo

\kern6pt
\line{\bf References\hfill}
\section{\,}\kern-20pt

\itemitem{[BL]}
{L. Biedenharn and J. Louck},
{A new class of symmetric
polynomials defined in terms of tableaux},
{Advances in Appl.\ Math.}
{\bf 10}
(1989),
396--438.

\itemitem{[C]}
{A. Capelli},
{Sur les op\'erations dans la th\'eorie des formes alg\'ebriques},
{Math. Ann.}
{\bf 37}
(1890),
1--37.

\itemitem{[C1]}
{I. Cherednik},
{Factorized particles on the half-line and root systems},
{Theor. Math. Phys.}
{\bf 61}
(1984),
977--983.

\itemitem{[C2]}
{I. Cherednik},
{On special bases of irreducible finite-dimensional representations
of the degenerate affine Hecke algebra},
{Funct. Anal. Appl.}
{\bf 20}
(1986),
87--89.

\itemitem{[D]}
{J. Dixmier},
{Alg\`ebres Enveloppantes}, 
{Gauthier-Villars, Paris},
1974.

\itemitem{[EHW]}
{T. Enright, R. Howe and N. Wallach},
{A classification of unitary highest weight modules},
in 
\lq\lq\tss Representation Theory of Reductive Groups\tss\rq\rq\,
(P. C. Trombi, Ed.),
Progress in Math. {\bf 40},
Birkh\"auser,
Bos\-ton,
1983,
pp. 97--143.

\itemitem{[H]}
{R. Howe},
{Remarks on classical invariant theory},
{Trans. Amer. Math. Soc.}
{\bf 313}
(1989),
539--570.

\itemitem{[HU]}
{R. Howe and T. Umeda},
{The Capelli identity, the double commutant theorem,
and multiplicity-free actions},
{Math. Ann.}
{\bf 290}
(1991),
569--619.

\itemitem{[K]}
{G. Kuperberg},
{Symmetries of plane partitions and the permanent-determinant  
method},
J. Comb. Theory
{\bf A\,68} 
(1994),
115--151.

\itemitem{[K1]}
F. Knop,
{A Harish-Chandra homomorphism for reductive group actions}, 
Ann. Math. 
{\bf 140} 
(1994), 
253--288.

\itemitem{[K2]}
F. Knop,
{Symmetric and non-symmetric quantum Capelli polynomials}, 
Comment. Math. Helvet.
{\bf 72}
(1997),
84--100.

\itemitem{[KS1]}
B. Kostant and S. Sahi,
{The Capelli identity, tube domains and the generalized Laplace 
transform}, 
Adv. Math. 
{\bf 87} 
(1991), 
71--92.

\itemitem{[KS2]}
B. Kostant and S. Sahi,
{Jordan algebras and Capelli identities},
Invent. Math. 
{\bf 112} 
(1993), 
657--664.

\itemitem{[KV]}
{M. Kashiwara and M. Vergne},
{On the Segal--Shale--Weil representations and harmonic polynomials},
Invent. Math.
{\bf 41}
(1978),
1--47.

\itemitem{[M]}
{I. Macdonald},
{Symmetric Functions and Hall Polynomials},
Clarendon Press, Oxford, 1995.

\itemitem{[M1]}
{A. Molev},
{Sklyanin determinant, Laplace operators, and characteristic identities
for classical Lie algebras}, 
{J. Math. Phys.}
{\bf 36}
(1995),
923--943.

\itemitem{[M2]}
{A. Molev},
{Factorial supersymmetric Schur functions and super Capelli identities},
in
\lq\lq\tss Kirillov's Seminar on Representation Theory\rq\rq\,
(G.~I.~Olshanski, Ed.),
{Amer. Math. Soc. Transl.}
{\bf 181},
AMS, Providence, 1998, pp. 109--137.

\itemitem{[MNO]}
{A. Molev, M. Nazarov and G. Olshanski},
{Yangians and classical Lie algebras}, 
Russian Math. Surveys
{\bf 51}
(1996),
205--282.

\itemitem{[N1]}
{M. Nazarov},
{Quantum Berezinian and the classical Capelli identity}, 
{Lett. Math. Phys.}
{\bf 21}
(1991),
123--131.

\itemitem{[N2]}
{M. Nazarov},
{Yangians and Capelli identities},
in
\lq\lq\tss Kirillov's Seminar on Representation Theory\rq\rq\,
(G.~I.~Olshanski, Ed.),
{Amer. Math. Soc. Transl.}
{\bf 181},
AMS, Providence, 1998, pp. 139--164.

\itemitem{[NO]}
{M. Nazarov and G. Olshanski},
{Bethe subalgebras in twisted Yangians}, 
Commun. Math. Phys.
{\bf 178}
(1996),
483--506.

\itemitem{[NUW1]}
{M. Noumi, T. Umeda and M. Wakayama},
{A quantum analogue of the Capelli identity
and an elementary differential calculus on $GL_q(n)$}, 
{Duke. Math. J.}
{\bf 76}
(1994),
567--594.

\itemitem{[NUW2]}
{M. Noumi, T. Umeda and M. Wakayama},
{Dual pairs, \text{spherical harmonics and} a Capelli identity
in quantum group theory}, 
{Compositio Math.}
{\bf 104}
(1996),
227--277.

\itemitem{[O]}
{G. Olshanski},
{Twisted Yangians and infinite-dimensional classical Lie algebras},
in 
\lq\lq\tss Quantum Groups\tss\rq\rq\,
(P. P. Kulish, Ed.),
{Lecture Notes in Math.} {\bf 1510},
Springer,
Berlin-Heidelberg,
1992, 
pp. 103--120.

\itemitem{[O1]}
{A. Okounkov},
{Quantum immanants and higher Capelli identities},
{Transformation Groups}
{\bf 1}
(1996),
99--126.

\itemitem{[O2]}
{A. Okounkov},
{Young basis, Wick formula, and higher Capelli identities},
{Int. Math. Research Notes}
(1996), 
817--839.

\itemitem{[OO1]}
A. Okounkov and G. Olshanski,
{Shifted Schur functions},
St.\,Petersburg Math. J.
{\bf 9}
(1998), 239--300.

\itemitem{[OO2]}
A. Okounkov and G. Olshanski,
{Shifted Schur functions II. 
Binomial formula for characters of classical groups and applications},
in
\lq\lq\tss Kirillov's Seminar on Representation Theory\rq\rq\,
(G.~I.~Olshanski, Ed.),
{Amer. Math. Soc. Transl.}
{\bf 181},
AMS, Providence, 1998, pp. 245--271.

\itemitem{[S]}
{S. Sahi},
{The spectrum of certain invariant differential operators associated 
to a Hermitian symmetric space,}
in
\lq\lq\tss Lie Theory and Geometry\tss\rq\rq\,
(J.-L. Brylinski, R. Brylinski, V. Guillemin, V. Kac, Eds.), 
{Progress in Math.} {\bf 123},
Birkh\"auser,
Boston,
1994,
pp. 569--576.

\itemitem{[S1]}
E. Sklyanin,
{Boundary conditions for integrable equations,} 
{Funct. Anal. Appl.}
{\bf 21}
(1987),
164--166.

\itemitem{[S2]}
E. Sklyanin,
{Boundary conditions for integrable quantum systems,} 
J. Phys.
{\bf A\,21}
(1988)
2375--2389.

\itemitem{[W]}
{H. Weyl},
{Classical Groups, their Invariants and Representations},
Princeton University Press, Princeton, 1946.


\bye